# MODERATE DEVIATIONS AND LAW OF THE ITERATED LOGARITHM FOR INTERSECTIONS OF THE RANGES OF RANDOM WALKS[1]

BY XIA CHEN

*University of Tennessee*

Let $S_1(n),\dots,S_p(n)$ be independent symmetric random walks in $\mathbb{Z}^d$. We establish moderate deviations and law of the iterated logarithm for the intersection of the ranges

$$\#\{S_1[0,n]\cap\cdots\cap S_p[0,n]\}$$

in the case $d=2$, $p\geq 2$ and the case $d=3$, $p=2$.

## 1. Introduction.

Let $p\geq 2$ be an integer and let $\{S_1(n)\},\dots,\{S_p(n)\}$ be symmetric independent $d$-dimensional lattice valued random walks with the same distribution. Throughout we assume that $\{S_1(n)\},\dots,\{S_p(n)\}$ have finite second moment and that the smallest group that supports these random walks is $\mathbb{Z}^d$. Write $\Gamma$ for their covariance matrix. Unless claiming otherwise, we assume that the random walks start at the origin, that is,

$$S_j(0)=0,\qquad j=1,\dots,p.$$

To simplify the notation, we use $\{S(n)\}$ for a random walk of the same distribution as $\{S_1(n)\},\dots,\{S_p(n)\}$, in the context where only a single random walk is involved. For any $\Delta\in\mathbb{R}^+$, we set

$$S(\Delta)=\{S(k);\ k\in\Delta\}.$$

In the transient case $d\geq 3$, we write

$$\gamma(S)=\mathbb{P}\{S(n)\neq 0,\ n\geq 1\}.$$

It is known [Dvoretzky, Erdös and Kakutani (1950, 1954)] that the trajectories of the random walks $\{S_1(n)\},\dots,\{S_p(n)\}$ intersect infinitely often

Received January 2004; revised May 2004.
[1]Supported in part by NSF Grant DMS-04-05188.

if and only if $p(d-2) \leq d$. There are two ways to measure the intensity of such intersection. One is to count the times of intersection by introducing the intersection local time

$$I_n = \#\{(k_1, \ldots, k_p) \in [0,n]^p;\ S_1(k_1) = \cdots = S_p(k_p)\}. \tag{1.1}$$

Another is to count the sites of intersection by considering the intersection of the ranges

$$J_n = \#\{S_1[0,n] \cap \cdots \cap S_p[0,n]\}. \tag{1.2}$$

In the critical cases defined by $p(d-2) = d$, a weak law obtained by Le Gall (1986b) shows that $I_n$ and $J_n$ are attracted by $\Gamma$-distributions. The law of the iterated logarithm (LIL) for $I_n$ and $J_n$ has been obtained in Marcus and Rosen (1997) and Rosen (1997). See (1.19) and (1.20) below for the LIL for $J_n$.

In Chen and Li (2004) and Chen (2004), the moderate deviations and the law of the iterated logarithm for $I_n$ have been established in the non-critical cases defined by $p(d-2) < d$. See also Chen, Li and Rosen (2005) and Chen and Rosen (2005) for the extensions of such results to the stable random walks.

In this paper, we study the moderate deviations and the law of the iterated logarithm for $J_n$ under the condition

$$p(d-2) < d \quad \text{and} \quad d \geq 2 \tag{1.3}$$

which consists of the case $d = 2$, $p \geq 2$ and the case $d = 3$, $p = 2$. Our work is partially inspired by two papers. One is Le Gall (1986a) in which it is pointed out [Theorem 5.1, Le Gall (1986b)] that as $d = 2$, $p \geq 2$, $m = 1, 2, \ldots,$

$$\frac{(\log n)^{pm}}{n^m} \mathbb{E} J_n^m \longrightarrow (2\pi)^{pm} \det(\Gamma)^{m/2} \mathbb{E}\alpha([0,1]^p)^m \qquad (n \to \infty) \tag{1.4}$$

and [Theorem 5.3, Le Gall (1986a)] that as $d = 3$ and $p = 2$, $m = 1, 2, \ldots,$

$$n^{-m/2} \mathbb{E} J_n^m \longrightarrow \gamma(S)^{2m} \det(\Gamma)^{-m/2} \mathbb{E}\alpha([0,1]^2)^m \qquad (n \to \infty) \tag{1.5}$$

where $\alpha([0,1]^p)$ is the Brownian intersection local time

$$\alpha([0,1]^p) = \int_{\mathbb{R}^d} \left[\prod_{j=1}^p \int_0^1 \delta_x(W_j(s))\, ds\right] dx \tag{1.6}$$

generated by the independent $d$-dimensional Brownian motions $W_1(t), \ldots, W_p(t)$. Here we make the following remarks: First, Le Gall only discussed the case where the covariance matrix $\Gamma$ is a multiple of the identical matrix. By examining his argument, we made a slight extension without repeating his proof. Second, it is very likely that (1.4) and (1.5) can be developed into the laws of weak convergence. To our best knowledge, this was confirmed [see, e.g., Le Gall (1986a) and Le Gall and Rosen (1991)] in the case $d=2$, $p=2,3$ and the case $d=3$ and $p=2$.

Another is the recent large deviation result [Theorem 2.1, Chen (2004); see also Chen and Rosen (2005) for its stable extension]

$$\lim_{t\to\infty} t^{-1}\log\mathbb{P}\{\alpha([0,1]^p)\geq t^{d(p-1)/2}\} = -\frac{p}{2}\kappa(d,p)^{-4p/(d(p-1))} \tag{1.7}$$

under the condition (1.3), where $\kappa(d,p)>0$ is the Gagliardo–Nirenberg constant given below. In view of (1.4) and (1.5), it is natural to expect that the tail behavior given in (1.7) passes to $J_n$ in certain ways.

For each $d$, $p$ satisfying (1.3), we introduce the positive number $\kappa(d,p)$ as the best constant of the Gagliardo–Nirenberg inequality

$$\|f\|_{2p}\leq C\|\nabla f\|_2^{d(p-1)/(2p)}\cdot\|f\|_2^{1-d(p-1)/(2p)},\qquad f\in W^{1,2}(\mathbb{R}^d),$$

where $W^{1,2}(\mathbb{R}^d)$ denotes the Sobolev space

$$W^{1,2}(\mathbb{R}^d)=\{f\in\mathcal{L}^2(\mathbb{R}^d);\ \nabla f\in\mathcal{L}^2(\mathbb{R}^d)\}.$$

That is,

$$\begin{aligned}\kappa(d,p)=\inf\{C>0;\ \|f\|_{2p}\leq C\|\nabla f\|_2^{d(p-1)/(2p)}&\cdot\|f\|_2^{1-d(p-1)/(2p)}\\ &\text{for } f\in W^{1,2}(\mathbb{R}^d)\}.\end{aligned} \tag{1.8}$$

The Gagliardo–Nirenberg inequality can be obtained from the Sobolev inequality by a simple substitution. We refer the interested reader to Levine (1980), Weinstein (1983), Carlen and Loss (1993), Del Pino and Dolbeault (2003) and Cordero-Erausquin, Nazaret and Villani (2004) for an overview of the latest state in finding the value of Gagliardo–Nirenberg constants.

THEOREM 1. *As $d=2$ and $p\geq 2$,*

$$\begin{aligned}\lim_{n\to\infty}\frac{1}{b_n}\log\mathbb{P}\Big\{J_n\geq\lambda\frac{n}{(\log n)^p}b_n^{p-1}\Big\}\\ =-\frac{p}{2}(2\pi)^{-p/(p-1)}\\ \times\det(\Gamma)^{-1/(2(p-1))}\kappa(2,p)^{-2p/(p-1)}\lambda^{1/(p-1)}\qquad(\lambda>0)\end{aligned} \tag{1.9}$$

*for each positive sequence $\{b_n\}$ satisfying*

$$b_n\to\infty\quad\text{and}\quad b_n=o((\log n)^{2/3})\qquad(n\to\infty). \tag{1.10}$$

THEOREM 2. *As $d=3$ and $p=2$,*

$$\begin{aligned}\lim_{n\to\infty}\frac{1}{b_n}\log\mathbb{P}\{J_n\geq\lambda\sqrt{nb_n^3}\}\\ =-\det(\Gamma)^{1/3}\gamma(S)^{-4/3}\kappa(3,2)^{-8/3}\lambda^{2/3}\qquad(\lambda>0)\end{aligned} \tag{1.11}$$

*for each positive sequence* $\{b_n\}$ *satisfying*

$$b_n \to \infty \quad \textit{and} \quad b_n = o(n^{2/9}) \qquad (n \to \infty). \tag{1.12}$$

REMARK. We point out the fact that as $d \geq 3$,

$$\gamma(S) = \left( \sum_{k=0}^{\infty} \mathbb{P}\{S(k) = 0\} \right)^{-1} = \left( \frac{1}{(2\pi)^d} \int_{[-\pi,\pi]^d} \frac{1}{1 - \varphi(\lambda)}\, d\lambda \right)^{-1} \tag{1.13}$$

where $\varphi(\lambda)$ is the characteristic function of the i.i.d. increments of $\{S(n)\}$. To prove the first equality in (1.13), let $\tau_0$ be the last time that the random walk $S(n)$ visits 0. By transience and the Markov property,

$$1 = \sum_{k=0}^{\infty} \mathbb{P}\{\tau_0 = k\} = \sum_{k=0}^{\infty} \mathbb{P}\{S(k) = 0\}\gamma(S).$$

The second equality in (1.13) follows from the fact that

$$\mathbb{P}\{S(k) = 0\} = \frac{1}{(2\pi)^d} \int_{[-\pi,\pi]^d} \varphi(\lambda)^k\, d\lambda, \qquad k = 0, 1, \ldots.$$

We now compare $J_n$ with $I_n$. A trivial observation gives that $J_n \leq I_n$ with the difference caused by the possibility that the multiple intersection may happen at the same site. By Theorem 2.2 in Chen (2004),

$$\lim_{n\to\infty} \frac{1}{b_n} \log \mathbb{P}\{I_n \geq \lambda n b_n^{p-1}\} = -\frac{p}{2} \sqrt{\det(\Gamma)} \kappa(2,p)^{-2p/(p-1)} \lambda^{(p-1)^{-1}} \tag{1.14}$$

as $d = 2$, $p \geq 2$; and

$$\lim_{n\to\infty} \frac{1}{b_n} \log \mathbb{P}\{I_n \geq \lambda \sqrt{n b_n^3}\} = -\det(\Gamma)^{1/3} \kappa(3,2)^{-8/3} \lambda^{2/3} \tag{1.15}$$

as $d = 3$, $p = 2$, where $\{b_n\}$ can be any positive sequence satisfying

$$b_n \to \infty \quad \text{and} \quad b_n = o(n) \qquad (n \to \infty). \tag{1.16}$$

Comparing (1.9) with (1.14), we see a substantial difference in asymptotic behaviors between $I_n$ and $J_n$ as $d = 2$.

Another difference is in the range of $\{b_n\}$. By comparison it is natural to ask if we can extend Theorems 1 and 2 so that any sequence $\{b_n\}$ satisfying (1.16) can be included. The answer is "No." Indeed, if we take $b_n \geq \delta(\log n)^{p/(p-1)}$ in (1.9), or $b_n \geq \delta n^{1/3}$ in (1.11), then the involved probability is bounded by

$$\mathbb{P}\{J_n \geq \delta\lambda n\}$$

which is eventually zero for $\lambda > \delta^{-1}$. So our results do not hold in this case.

It seems that in Theorem 2, the right condition on $\{b_n\}$ is

$$b_n \to \infty \quad \text{and} \quad b_n = o(n^{1/3}) \qquad (n \to \infty).$$

As for Theorem 1, we can push a little further: If (1.9) were true for $b_n = \log n$, we would have

$$\lim_{n\to\infty} \frac{1}{\log n} \log \mathbb{P}\left\{ J_n \geq \lambda \frac{n}{\log n} \right\}$$
$$= -\frac{p}{2}(2\pi)^{-p/(p-1)} \det(\Gamma)^{-1/(2(p-1))} \kappa(2,p)^{-2p/(p-1)} \lambda^{1/(p-1)}.$$

This is implausible since, in the sense of moderate deviation at the scale $b_n = \log n$, $J_n$ would have the rate $n(\log n)^{-1}$ independent of $p$, which sharply contrasts with (1.4). We believe that in Theorem 1, the right condition on $\{b_n\}$ is

$$b_n \to \infty \quad \text{and} \quad b_n = o(\log n) \qquad (n \to \infty).$$

We are not able to prove our results under these conditions. So we leave this problem to future study.

THEOREM 3. *As $d = 2$ and $p \geq 2$,*

$$\limsup_{n\to\infty} \frac{(\log n)^p}{n(\log\log n)^{p-1}} J_n = (2\pi)^p \left(\frac{2}{p}\right)^{p-1} \sqrt{\det(\Gamma)}\kappa(2,p)^{2p} \qquad a.s. \tag{1.17}$$

*As $d = 3$ and $p = 2$,*

$$\limsup_{n\to\infty} \frac{1}{\sqrt{n(\log\log n)^3}} J_n = \gamma(S)^2 \det(\Gamma)^{-1/2} \kappa(3,2)^4 \qquad a.s. \tag{1.18}$$

Recall that the trajectories of $\{S_1(n)\}, \ldots, \{S_p(n)\}$ intersect infinitely often if and only if $p(d-2) \leq d$. In the critical cases defined as $p(d-2) = d$—the case "$d = 4$, $p = 2$" and the case "$d = p = 3$," the law of the iterated logarithm for $J_n$ has been obtained in Marcus and Rosen (1997) and in Rosen (1997), respectively. Under the assumption of finite third moment, it has been proved [Marcus and Rosen (1997)] that

$$\limsup_{n\to\infty} \frac{J_n}{\log n \log\log\log n} = \frac{\gamma(S)^2}{2\pi^2\sqrt{\det(\Gamma)}} \qquad \text{a.s.} \tag{1.19}$$

as $d = 4$ and $p = 2$, and [Rosen (1997)] that

$$\limsup_{n\to\infty} \frac{J_n}{\log n \log\log\log n} = \frac{\gamma(S)^3}{\pi\det(\Gamma)} \qquad \text{a.s.} \tag{1.20}$$

as $d = p = 3$.

As $d = 1$, we have

$$J_n \le \min_{1\le j\le p} \max_{k\le n} S_j(k) - \max_{1\le j\le p} \min_{k\le n} S_j(k). \tag{1.21}$$

Since the equality holds in the special case of simple random walks, it is natural to believe that even in the general case, both sides of (1.21) are asymptotically equivalent in a suitable sense. By the classical results on the tail estimate of the random walks, therefore, we conjecture that

$$\lim_{n\to\infty} \frac{1}{b_n} \log \mathbb{P}\{J_n \ge \lambda\sqrt{nb_n}\} = -\frac{p\lambda^2}{2\sigma^2} \tag{1.22}$$

for any positive sequence $\{b_n\}$ satisfying (1.16), where $\sigma^2 > 0$ is the variance of the random walks. The rigorous proof of (1.22) [more precisely, the lower bound of (1.22)] for the general random walks can be difficult. By comparing (1.22) with Theorems 1, 2, (1.19) and (1.20), it is interesting to note that the asymptotic magnitude of $J_n$ is not monotonic in dimension $d$ and that asymptotically, $J_n$ is maximized by $d = 2$.

Another interesting problem is the study of $\#\{S[0,n]\}$ (i.e., $J_n$ with $p = 1$). In the case $d = 1$, it is expected that $\#\{S[0,n]\}$ behaves like

$$\max_{k\le n} S(k) - \min_{k\le n} S(k)$$

in terms of the upper and lower tail behaviors.

In the multidimensional case, the behaviors of the range $\#\{S[0,n]\}$ are generally different from what we observe in the present paper. In the case $d \ge 3$, it has been shown [Jain and Pruitt (1972) and Bass and Kumagai (2002)] that the centered sequence

$$\#\{S[0,n]\} - \mathbb{E}\#\{S[0,n]\}, \qquad n = 1, 2, \dots, \tag{1.23}$$

has Gaussian tails and behaves essentially like a partial sum of independent random variables.

The case $d = 2$ is the most interesting case in which the tail of the sequence in (1.23) is no longer Gaussian, not even symmetric. Bass and Kumagai (2002) obtain

$$\limsup_{n\to\infty} \frac{(\log n)^2}{n \log\log\log n} (\#\{S[0,n]\} - \mathbb{E}\#\{S[0,n]\}) = C \qquad \text{a.s.} \tag{1.24}$$

with the unidentified constant $C > 0$. In a forthcoming paper, we [Bass, Chen and Rosen (2004)] shall identify the constant $C$ and we shall show that it is the lim inf behavior of the sequence in (1.23) (i.e., $J_n - \mathbb{E}J_n$ with $p = 1$) that is relevant to the lim sup behavior of $J_n$ (with $p = 2$) given in Theorem 3.

Finally, we point out some interesting problems in the case

$$p(d-2) > d. \tag{1.25}$$

According to Dvoretzky, Erdös and Kakutani (1950, 1954), we have

$$I_\infty = \#\{(k_1,\dots,k_p)\in[0,\infty)^p;\ S_1(k_1)=\cdots=S_p(k_p)\}<\infty \qquad \text{a.s.},$$

$$J_\infty = \#\{S_1[0,\infty)\cap\cdots\cap S_p[0,\infty)\}<\infty \qquad \text{a.s.};$$

a natural problem is to study the tails of the random variables $I_\infty$ and $J_\infty$. In Khanin, Mazel, Shlosman and Sinai (1994), this problem is linked to the study of the random walk in the random potential. In the special case $d\ge 5$ and $p=2$, Khanin, Mazel, Shlosman and Sinai (1994) prove that there are $c_1, c_2>0$, such that

$$\exp\{-c_1 t^{1/2}\}\le \mathbb{P}\{I_\infty\ge t\}\le \exp\{-c_2 t^{1/2}\} \tag{1.26}$$

and that given $\delta>0$,

$$\exp\{-t^{1-2/d+\delta}\}\le \mathbb{P}\{J_\infty\ge t\}\le \exp\{-t^{1-2/d-\delta}\} \tag{1.27}$$

holds for large $t$. From (1.26) and (1.27) we observe again a fundamental difference between the intersection local time and the intersection of independent ranges. In particular, this observation breaks the stereotype that $J.$ always behaves like $\gamma(S)^p I.$ in the transient case. It is certainly of great interest in studying precise large deviations for $I_\infty$ and $J_\infty$ under (1.25).

The paper is organized as follows. In Section 2, we formulate a nonstandard version (Theorem 4) of the Gärtner–Ellis theorem with nearly standard proof. From the viewpoint of large deviation theory, our work contributes an important example which is not quite suitable for the classic Gärtner–Ellis theorem but can be solved in a nonstandard way.

In Section 3, we prove the upper bounds given in Theorems 1 and 2. The key tool is a moment inequality (Theorem 6) for $J_n$ which is parallel to the one given in Theorem 5.1 in Chen (2004) for $I_n$.

In Section 4, we prove the lower bounds given in Theorems 1 and 2. This is the most delicate part of the whole paper and some substantially new ideas are needed. First we establish a weak law (Theorem 7) for certain functionals related to $J_n$, which seems new and has independent interest for its own sake. Second, we partition the time interval $[0,n]$ properly and conduct some sharp estimate to eliminate the influence from intersection of trajectories between any two different time periods. Finally, we establish some Feynman–Kac type large deviation lower bounds (Theorem 8) in a way close to Theorem 4.1 in Chen and Li (2004).

In Section 5, we prove the laws of the iterated logarithm given in Theorem 3. The nontrivial part is the lower bound, for which some uniform lower bounds of the moderate deviations are needed.

In spite of some technical connections to the recent works Chen and Li (2004), Bass and Chen (2004), Chen (2004), Chen, Li and Rosen (2005), Chen and Rosen (2005) and Bass, Chen and Rosen (2005) on the exponential asymptotics for intersection local times, the main approach used here is fundamentally different.

**2. A Gärtner–Ellis type theorem.** Let $\{Z_\varepsilon\}$ be a family of nonnegative random variables and let $p \ge 1$ be an integer. Assume that for any $\theta > 0$, the following limit exists:

$$\lim_{\varepsilon \to 0^+} \varepsilon \log \sum_{m=0}^{\infty} \frac{(\theta \varepsilon^{-1})^m}{m!} (\mathbb{E} Z_\varepsilon^m)^{1/p} = \Psi(\theta). \tag{2.1}$$

It is easy to see that $\Psi(\theta)$ is nondecreasing and convex on $[0, \infty)$ with $\Psi(0) = 0$. By the Gärtner–Ellis theorem, $Z_\varepsilon$ satisfies the large deviation principle if $p = 1$ and if $\Psi(\theta)$ and its convex conjugate $\Psi^*(\lambda)$ satisfy some regularity conditions [see, e.g., Theorem 2.3.6 in Dembo and Zeitouni (1998) for details]. What we intend to establish in this section is a large deviation principle under (2.1) and some additional regularity assumptions in the case $p \ge 1$.

Write

$$I(\lambda) = p \sup_{\theta > 0} \{\lambda^{1/p} \theta - \Psi(\theta)\}. \tag{2.2}$$

By Lemma 2.3.9 in Dembo and Zeitouni (1998), $I$ is a good rate function: $I$ is lower semicontinuous on $[0, \infty]$ and for each $l > 0$, the level set $\{\lambda; I(\lambda) \le l\}$ is compact. In addition, one can easily see that $I(0) = 0$ and that $I(|\cdot|^p)$ is convex on $(-\infty, \infty)$.

Definition. $\lambda_0 \in [0, \infty)$ is called a $p$-distinguishable point of $I$ if there is $\theta_0 \in [0, \infty)$ such that

$$\lambda^{1/p} \theta_0 - \frac{1}{p} I(\lambda) < \Psi(\theta_0) \qquad \forall \lambda > 0 \text{ with } \lambda \ne \lambda_0.$$

Remark. By an argument of duality [see the proof of Lemma 5.3 in Chen (2004)] we have that for any $\theta_0 > 0$,

$$\sup_{\lambda > 0} \left\{ \lambda^{1/p} \theta_0 - \frac{1}{p} I(\lambda) \right\} = \Psi(\theta_0).$$

Therefore, $\lambda_0$ is $p$-distinguishable if $\lambda_0$ is the unique maximizer of the function

$$\varphi(\lambda) = \lambda^{1/p} \theta_0 - \frac{1}{p} I(\lambda)$$

for some $\theta_0 \ge 0$.

An important ingredient of our idea is the following generalization of the Gärtner–Ellis theorem on large deviations.

THEOREM 4. *Let $\{Z_\varepsilon\}$ be a family of nonnegative random variables and let $p \ge 1$ be an integer. Assume that for any $\theta > 0$, (2.1) holds. Then for any $\lambda > 0$,*

$$\limsup_{\varepsilon \to 0^+} \varepsilon \log \mathbb{P}\{Z_\varepsilon \ge \lambda\} \le -I(\lambda). \tag{2.3}$$

*Further, if the set of $p$-distinguishable points of $I$ is dense in $[0, \infty)$, then*

$$\lim_{\varepsilon \to 0^+} \varepsilon \log \mathbb{P}\{Z_\varepsilon \ge \lambda\} = -I(\lambda), \qquad \lambda > 0. \tag{2.4}$$

PROOF. The proof of the upper bound is just a routine application of the Chebyshev inequality: For any $\theta > 0$,

$$\lambda^{m/p}(\theta\varepsilon^{-1})^m (\mathbb{P}\{Z_\varepsilon \ge \lambda\})^{1/p} \le (\theta\varepsilon^{-1})^m (\mathbb{E}Z_\varepsilon^m)^{1/p}$$

for any integer $m \ge 0$. Summing up gives

$$e^{\theta\lambda^{1/p}\varepsilon^{-1}}(\mathbb{P}\{Z_\varepsilon \ge \lambda\})^{1/p} \le \sum_{m=0}^{\infty} \frac{(\theta\varepsilon^{-1})^m}{m!}(\mathbb{E}Z_\varepsilon^m)^{1/p}.$$

Hence

$$\limsup_{\varepsilon \to 0^+} \varepsilon \log \mathbb{P}\{Z_\varepsilon \ge \lambda\} \le -p\{\lambda^{1/p}\theta - \Psi(\theta)\}.$$

Taking the supremum over $\theta$ gives the desired upper bound.

To accomplish the second part, we need only to prove that for any $p$-distinguishable point $\lambda_0$ and any $\delta > 0$,

$$\liminf_{\varepsilon \to 0^+} \varepsilon \log \mathbb{P}\{Z_\varepsilon \in (\lambda_0 - \delta, \lambda_0 + \delta)\} \ge -I(\lambda_0). \tag{2.5}$$

We may assume that $0 < \delta < \lambda_0$. Notice that

$$(\lambda_0 + \delta)^{m/p}(\mathbb{P}\{Z_\varepsilon \in (\lambda_0 - \delta, \lambda_0 + \delta)\})^{1/p} \ge (\mathbb{E}Z_\varepsilon^m \mathbb{1}_{\{Z_\varepsilon \in (\lambda_0 - \delta, \lambda_0 + \delta)\}})^{1/p}.$$

Summing up we have

$$\begin{aligned} &e^{\varepsilon^{-1}\theta_0(\lambda_0+\delta)^{1/p}}(\mathbb{P}\{Z_\varepsilon \in (\lambda_0 - \delta, \lambda_0 + \delta)\})^{1/p} \\ &\quad \ge \sum_{m=0}^{\infty} \frac{(\theta_0\varepsilon^{-1})^m}{m!}(\mathbb{E}Z_\varepsilon^m \mathbb{1}_{\{Z_\varepsilon \in (\lambda_0 - \delta, \lambda_0 + \delta)\}})^{1/p} \end{aligned}$$

where $\theta_0$ is given as in the definition of the $p$-distinguishable point $\lambda_0$.

If we can prove that for any $\delta > 0$,

$$\begin{aligned} &\sum_{m=0}^{\infty} \frac{(\theta_0\varepsilon^{-1})^m}{m!}(\mathbb{E}Z_\varepsilon^m \mathbb{1}_{\{Z_\varepsilon \in (\lambda_0 - \delta, \lambda_0 + \delta)\}})^{1/p} \\ &\quad \sim \sum_{m=0}^{\infty} \frac{(\theta_0\varepsilon^{-1})^m}{m!}(\mathbb{E}Z_\varepsilon^m)^{1/p} \qquad (\varepsilon \to 0^+) \end{aligned} \tag{2.6}$$

then we will have

$$\liminf_{\varepsilon\to 0^+} \varepsilon \log \mathbb{P}\{Z_\varepsilon \in (\lambda_0-\delta, \lambda_0+\delta)\} \ge -p\{\theta_0(\lambda+\delta)^{1/p} - \Psi(\theta_0)\}.$$

For any $0<\delta'<\delta$, replacing $\delta$ by $\delta'$ and noticing that

$$\mathbb{P}\{Z_\varepsilon \in (\lambda_0-\delta, \lambda_0+\delta)\} \ge \mathbb{P}\{Z_\varepsilon \in (\lambda_0-\delta', \lambda_0+\delta')\},$$

we obtain

$$\liminf_{\varepsilon\to 0^+} \varepsilon \log \mathbb{P}\{Z_\varepsilon \in (\lambda_0-\delta, \lambda_0+\delta)\} \ge -p\{\theta_0(\lambda+\delta')^{1/p} - \Psi(\theta_0)\}.$$

Letting $\delta' \to 0^+$ gives

$$\liminf_{\varepsilon\to 0^+} \varepsilon \log \mathbb{P}\{Z_\varepsilon \in (\lambda_0-\delta, \lambda_0+\delta)\} \ge -p\{\theta_0\lambda^{1/p} - \Psi(\theta_0)\} \ge -I(\lambda_0).$$

That is (2.5).

To prove (2.6), notice that

$$\begin{aligned}\sum_{m=0}^\infty &\frac{(\theta_0\varepsilon^{-1})^m}{m!}(\mathbb{E}Z_\varepsilon^m)^{1/p}\\ &\le \sum_{m=0}^\infty \frac{(\theta_0\varepsilon^{-1})^m}{m!}(\mathbb{E}Z_\varepsilon^m \mathbb{1}_{\{Z_\varepsilon\in(\lambda_0-\delta,\lambda_0+\delta)\}})^{1/p}\\ &\quad + \sum_{m=0}^\infty \frac{(\theta_0\varepsilon^{-1})^m}{m!}(\mathbb{E}Z_\varepsilon^m \mathbb{1}_{\{Z_\varepsilon\notin(\lambda_0-\delta,\lambda_0+\delta)\}})^{1/p}.\end{aligned}$$

In view of (2.1), we will have (2.6) if

$$\liminf_{\varepsilon\to 0^+} \varepsilon \log \sum_{m=0}^\infty \frac{(\theta_0\varepsilon^{-1})^m}{m!}(\mathbb{E}Z_\varepsilon^m \mathbb{1}_{\{Z_\varepsilon\notin(\lambda_0-\delta,\lambda_0+\delta)\}})^{1/p} < \Psi(\theta_0). \tag{2.7}$$

Write $B_0 = (\lambda_0-\delta, \lambda_0+\delta)$. Since $I(\lambda)$ is a good rate function, by distinguishability

$$\eta \equiv \Psi(\theta_0) - \sup_{\lambda\notin B_0}\left\{\lambda^{1/p}\theta_0 - \frac{1}{p}I(\lambda)\right\} > 0.$$

From the Hölder inequality, $(\mathbb{E}Z_\varepsilon^m)^{1/p} \ge \mathbb{E}Z_\varepsilon^{m/p}$ and the assumption (2.1) we have

$$\limsup_{\varepsilon\to 0} \varepsilon \log \mathbb{E}\exp\{\theta\varepsilon^{-1}Z_\varepsilon\} < \infty, \qquad \theta > 0.$$

According to Lemma 5.3(iii) in Chen (2004) (or Theorem 5 below), therefore,

$$\lim_{N\to\infty}\limsup_{\varepsilon\to 0^+} \varepsilon \log \sum_{m=0}^\infty \frac{(\theta_0\varepsilon^{-1})^m}{m!}(\mathbb{E}Z_\varepsilon^m \mathbb{1}_{\{Z_\varepsilon\ge N\}})^{1/p} = -\infty.$$

Let $N > \lambda + \delta$ be fixed for a moment and let

$$B_i = [a_i, b_i], \qquad i = 1, \dots, l,$$

be intervals such that

$$[0, N] \setminus B_0 = \bigcup_{i=1}^{l} B_i$$

and that $(b_i - a_i)^{1/p} < \eta/2$ for $i = 1, \dots, l$. Then

$$\begin{aligned}
&\sum_{m=0}^{\infty} \frac{(\theta_0 \varepsilon^{-1})^m}{m!} (\mathbb{E} Z_\varepsilon^m \mathbb{1}_{\{Z_\varepsilon \notin B_0\}})^{1/p} \\
&\quad \le \sum_{m=0}^{\infty} \frac{(\theta_0 \varepsilon^{-1})^m}{m!} (\mathbb{E} Z_\varepsilon^m \mathbb{1}_{\{Z_\varepsilon \ge N\}})^{1/p} + \sum_{i=1}^{l} \sum_{m=0}^{\infty} \frac{(\theta_0 \varepsilon^{-1})^m}{m!} (\mathbb{E} Z_\varepsilon^m \mathbb{1}_{\{Z_\varepsilon \in B_i\}})^{1/p} \\
&\quad \le \sum_{m=0}^{\infty} \frac{(\theta_0 \varepsilon^{-1})^m}{m!} (\mathbb{E} Z_\varepsilon^m \mathbb{1}_{\{Z_\varepsilon \ge N\}})^{1/p} + \sum_{i=1}^{l} e^{\theta_0 b_i \varepsilon^{-1}} (\mathbb{P}\{Z_\varepsilon \ge a_i\})^{1/p}.
\end{aligned}$$

By the proved upper bound,

$$\begin{aligned}
&\limsup_{\varepsilon \to 0^+} \varepsilon \log \sum_{m=0}^{\infty} \frac{(\theta_0 \varepsilon^{-1})^m}{m!} (\mathbb{E} Z_\varepsilon^m \mathbb{1}_{\{Z_\varepsilon \notin B_0\}})^{1/p} \\
&\quad \le \max\left\{ \limsup_{\varepsilon \to 0^+} \varepsilon \sum_{m=0}^{\infty} \frac{(\theta_0 \varepsilon^{-1})^m}{m!} (\mathbb{E} Z_\varepsilon^m \mathbb{1}_{\{Z_\varepsilon \ge N\}})^{1/p}, \right. \\
&\qquad\qquad \left. \max_{1 \le i \le l} \left\{ \theta_0 b_i^{1/p} - \frac{1}{p} I(a_i) \right\} \right\} \\
&\quad \le \max\left\{ \limsup_{\varepsilon \to 0^+} \varepsilon \sum_{m=0}^{\infty} \frac{(\theta_0 \varepsilon^{-1})^m}{m!} (\mathbb{E} Z_\varepsilon^m \mathbb{1}_{\{Z_\varepsilon \ge N\}})^{1/p}, \right. \\
&\qquad\qquad \left. \sup_{\lambda \notin B_0} \left\{ \theta_0 \lambda^{1/p} - \frac{1}{p} I(\lambda) \right\} + \frac{\eta}{2} \right\}.
\end{aligned}$$

Letting $N \to \infty$ gives

$$\begin{aligned}
&\limsup_{\varepsilon \to 0^+} \varepsilon \log \sum_{m=0}^{\infty} \frac{(\theta_0 \varepsilon^{-1})^m}{m!} (\mathbb{E} Z_\varepsilon^m \mathbb{1}_{\{Z_\varepsilon \notin B_0\}})^{1/p} \\
&\quad \le \sup_{\lambda \notin B_0} \left\{ \theta_0 \lambda^{1/p} - \frac{1}{p} I(\lambda) \right\} + \frac{\eta}{2} < \Psi(\theta_0).
\end{aligned}$$

□

Like Varadhan's integral lemma [Theorem 4.3.1 in Dembo and Zeitouni (1998)] to the well-known Gärtner–Ellis theorem, the following theorem is a converse of Theorem 4. We give it without proof, as it is essentially given in the proof for Lemma 5.3 in Chen (2004) (only some obvious modification is needed).

THEOREM 5. *Let* $\{Z_\varepsilon\}$ *be a family of nonnegative random variables and let* $p \geq 1$ *be an integer. Let* $I(\lambda)$ *be a nondecreasing good rate function on* $[0,\infty)$ *such that* $I(0)=0$, $I(|\cdot|^p)$ *is convex on* $(-\infty,\infty)$. *Assume that*

$$\lim_{\varepsilon\to 0^+} \varepsilon \log \mathbb{P}\{Z_\varepsilon \geq \lambda\} = -I(\lambda) \qquad (\lambda > 0)$$

*and that* $\theta > 0$ *satisfies*

$$\lim_{N\to\infty} \limsup_{\varepsilon\to 0^+} \varepsilon \log \sum_{m=1}^{\infty} \frac{(\theta\varepsilon^{-1})^m}{m!} (\mathbb{E} Z_\varepsilon^m \mathbb{1}_{\{Z_\varepsilon \geq N\}})^{1/p} = -\infty; \tag{2.8}$$

*then*

$$\lim_{\varepsilon\to 0^+} \varepsilon \log \left( \sum_{m=0}^{\infty} \frac{(\theta\varepsilon^{-1})^m}{m!} (\mathbb{E} Z_\varepsilon^m)^{1/p} \right) = \sup_{\lambda>0} \{\theta\lambda^{1/p} - p^{-1} I(\lambda)\}. \tag{2.9}$$

*In particular, the condition* (2.8) *is satisfied if there is a* $\theta' > 2p\theta$ *such that*

$$\limsup_{\varepsilon\to 0^+} \varepsilon \log \mathbb{E} \exp\{\varepsilon^{-1}\theta' Z_\varepsilon^{1/p}\} < \infty. \tag{2.10}$$

Theorem 4 applies to the proof of Theorems 1 and 2 as follows.

CLAIM 1. *We will have Theorem* 1 *if*

$$\begin{aligned} &\lim_{n\to\infty} \frac{1}{b_n} \log \sum_{m=0}^{\infty} \frac{\theta^m}{m!} \left(\frac{b_n \log^p n}{n}\right)^{m/p} (\mathbb{E} J_n^m)^{1/p} \\ &\quad = \frac{1}{p} \left(\frac{2(p-1)}{p}\right)^{p-1} (2\pi\theta)^p \sqrt{\det(\Gamma)}\, \kappa(2,p)^{2p} \qquad (\theta > 0) \end{aligned} \tag{2.11}$$

*in the case* $d=2$, $p \geq 2$.

CLAIM 2. *We will have Theorem* 2 *if*

$$\begin{aligned} &\lim_{n\to\infty} \frac{1}{b_n} \log \sum_{m=0}^{\infty} \frac{\theta^m}{m!} \left(\frac{b_n}{n}\right)^{m/4} \sqrt{\mathbb{E} J_n^m} \\ &\quad = 2\left(\frac{3}{4}\right)^3 (\gamma(S)\theta)^4 \det(\Gamma)^{-1} \kappa(3,2)^8 \qquad (\theta > 0) \end{aligned} \tag{2.12}$$

*in the case* $d=3$, $p=2$.

Due to similarity we only show how Claim 1 follows from Theorem 4. First, the condition (2.1) is satisfied with

$$\Psi(\theta) = \frac{1}{p}\left(\frac{2(p-1)}{p}\right)^{p-1}(2\pi\theta)^p\sqrt{\det(\Gamma)}\kappa(2,p)^{2p}.$$

A simple calculus gives that

$$\begin{aligned} I(\lambda) &= p\sup_{\theta>0}\{\lambda^{1/p}\theta - \Psi(\theta)\} \\ &= \frac{p}{2}(2\pi)^{-p/(p-1)}\det(\Gamma)^{-1/(2(p-1))}\kappa(2,p)^{-2p/(p-1)}\lambda^{1/(p-1)}. \end{aligned}$$

Second, every $\lambda_0 > 0$ is $p$-distinguishable. Indeed, doing simple calculus again one can directly verify that for

$$\theta_0 = \frac{1}{2}\frac{p}{p-1}(2\pi)^{-p/(p-1)}\kappa(2,p)^{-2p/(p-1)}\lambda_0^{1/(p(p-1))},$$

$\lambda_0$ is the unique maximizer of the function

$$\varphi(\lambda) = \lambda^{1/p}\theta_0 - \frac{1}{p}I(\lambda).$$

**3. Upper bounds.** The main goal of this section is to prove that in the case $d=2$, $p\geq 2$,

$$\begin{aligned} &\limsup_{n\to\infty}\frac{1}{b_n}\log\sum_{m=0}^{\infty}\frac{\theta^m}{m!}\left(\frac{b_n\log^p n}{n}\right)^{m/p}(\mathbb{E}J_n^m)^{1/p} \\ &\quad\leq \frac{1}{p}\left(\frac{2(p-1)}{p}\right)^{p-1}(2\pi\theta)^p\sqrt{\det(\Gamma)}\kappa(2,p)^{2p} \qquad (\theta>0) \end{aligned} \tag{3.1}$$

for any $\{b_n\}$ satisfying (1.10); and that in the case $d=3$, $p=2$,

$$\begin{aligned} &\limsup_{n\to\infty}\frac{1}{b_n}\log\sum_{m=0}^{\infty}\frac{\theta^m}{m!}\left(\frac{b_n}{n}\right)^{m/4}\sqrt{\mathbb{E}J_n^m} \\ &\quad\leq 2\left(\frac{3}{4}\right)^3(\gamma(S)\theta)^4\det(\Gamma)^{-1}\kappa(3,2)^8 \qquad (\theta>0) \end{aligned} \tag{3.2}$$

for any $\{b_n\}$ satisfying (1.12).

To begin, we first consider $\{S_1(n)\},\dots,\{S_p(n)\}$ as any independent and identically distributed $\mathbb{Z}^d$-random walks. Let the integer $a\geq 2$ be fixed and let $n_1,\dots,n_a$ be positive integers, $n_0=0$. Write

$$\Delta_i = [n_0+\cdots+n_{i-1},\, n_0+\cdots+n_i], \qquad i=1,\dots,a,$$

$$A = \sum_x\prod_{j=1}^p\sum_{i=1}^a \mathbb{1}_{\{x\in S_j(\Delta_i)\}}.$$

Notice that

$$J_{n_1+\cdots+n_a} = \sum_x \prod_{j=1}^p \mathbb{1}_{\{x\in S_j[0,n_1+\cdots+n_a]\}} \le A.$$

For the needs of the upper bound, it is enough to control $J_{n_1+\cdots+n_a}$. In the proof of the lower bound, however, it is required to control the self-intersection between two different parts of a single trajectory, which is associated with $A$ (with $a, n_1, \dots, n_a$ being suitably chosen) in law. In addition, the hardest part of this work is to essentially show that $A$ and $J_{n_1+\cdots+n_a}$ are asymptotically equivalent as $a, n_1, \dots, n_a$ (all depend on $n$) are suitably chosen.

THEOREM 6. *For any integer $m \ge 1$,*

$$(\mathbb{E}A^m)^{1/p} \le \sum_{\substack{k_1+\cdots+k_m=m\\ k_1,\dots,k_m\ge 0}} \frac{m!}{k_1!\cdots k_a!}(\mathbb{E}J_{n_1}^{k_1})^{1/p}\cdots(\mathbb{E}J_{n_a}^{k_a})^{1/p}. \tag{3.3}$$

*Consequently, for any $\lambda > 0$,*

$$\sum_{m=0}^\infty \frac{\theta^m}{m!}(\mathbb{E}A^m)^{1/p} \le \prod_{i=1}^a \sum_{m=0}^\infty \frac{\theta^m}{m!}(\mathbb{E}J_{n_i}^m)^{1/p}. \tag{3.4}$$

PROOF.

$$\begin{aligned}(\mathbb{E}A^m)^{1/p} &= \left(\sum_{x_1,\dots,x_m}\left[\mathbb{E}\prod_{k=1}^m\sum_{i=1}^a \mathbb{1}_{\{x_k\in S(\Delta_i)\}}\right]^p\right)^{1/p}\\ &= \left(\sum_{x_1,\dots,x_m}\left[\sum_{i_1,\dots,i_m=1}^a \mathbb{E}(\mathbb{1}_{\{x_1\in S(\Delta_{i_1})\}}\cdots\mathbb{1}_{\{x_m\in S(\Delta_{i_m})\}})\right]^p\right)^{1/p}\\ &\le \sum_{i_1,\dots,i_m=1}^a\left(\sum_{x_1,\dots,x_m}[\mathbb{E}(\mathbb{1}_{\{x_1\in S(\Delta_{i_1})\}}\cdots\mathbb{1}_{\{x_m\in S(\Delta_{i_m})\}})]^p\right)^{1/p}.\end{aligned}$$

Given integers $i_1, \dots, i_m$ between 1 and $a$, let $k_1, \dots, k_a$ be the number of occurrences of $i. = 1, \dots, i. = a$, respectively. Then $k_1 + \cdots + k_a = m$. To prove (3.3), it suffices to show

$$\sum_{x_1\cdots x_m}[\mathbb{E}(\mathbb{1}_{\{x_1\in S(\Delta_{i_1})\}}\cdots\mathbb{1}_{\{x_m\in S(\Delta_{i_m})\}})]^p \le \mathbb{E}J_{n_1}^{k_1}\cdots\mathbb{E}J_{n_a}^{k_a}. \tag{3.5}$$

Without losing generality we may only consider the case when $k_1, \dots, k_a \ge 1$. Under the notation $\bar{x}_i = (x_1^i, \dots, x_{k_i}^i) \in (\mathbb{Z}^d)^{k_i}$, we set

$$\phi_i(\bar{x}_i) = \mathbb{E}\left(\prod_{l=1}^{k_i}\mathbb{1}_{\{x_l^i\in S[0,n_i]\}}\right).$$

It is easy to see that

$$\sum_{\bar{x}_i} \phi_i^p(\bar{x}_i) = \mathbb{E} J_{n_i}^{k_i}, \qquad i = 1, \dots, a.$$

Define

$$\bar{S}^i(k) = (\overbrace{S(k), \dots, S(k)}^{k_i}) \quad \text{and} \quad \bar{S}_j^i(k) = (\overbrace{S_j(k), \dots, S_j(k)}^{k_i}), \qquad k = 1, 2, \dots,$$

where $1 \le i \le a$ and $1 \le j \le p$. Then

$$\begin{aligned} &\sum_{x_1 \cdots x_m} [\mathbb{E}(\mathbb{1}_{\{x_1 \in S(\Delta_{i_1})\}} \cdots \mathbb{1}_{\{x_m \in S(\Delta_{i_m})\}})]^p \\ &\qquad = \sum_{\bar{x}_1} \cdots \sum_{\bar{x}_a} \left[ \mathbb{E} \prod_{i=1}^a \mathbb{1}_{\{x_1^i \in S(\Delta_i)\}} \cdots \mathbb{1}_{\{x_{k_i}^i \in S(\Delta_i)\}} \right]^p. \end{aligned}$$

Notice that

$$\begin{aligned} &\sum_{\bar{x}_a} \left[ \mathbb{E} \prod_{i=1}^a \mathbb{1}_{\{x_1^i \in S(\Delta_i)\}} \cdots \mathbb{1}_{\{x_{k_i}^i \in S(\Delta_i)\}} \right]^p \\ &\quad = \sum_{\bar{x}_a} \left[ \mathbb{E} \left\{ \left( \prod_{i=1}^{a-1} \mathbb{1}_{\{x_1^i \in S(\Delta_i)\}} \cdots \mathbb{1}_{\{x_{k_i}^i \in S(\Delta_i)\}} \right) \phi_a(\bar{x}_a - \bar{S}^a(n - n_a)) \right\} \right]^p \\ &\quad = \sum_{\bar{x}_a} \mathbb{E} \left\{ \prod_{j=1}^p \left( \prod_{i=1}^{a-1} \mathbb{1}_{\{x_1^i \in S_j(\Delta_i)\}} \cdots \mathbb{1}_{\{x_{k_i}^i \in S_j(\Delta_i)\}} \right) \phi_a(\bar{x}_a - \bar{S}_j^a(n - n_a)) \right\} \\ &\quad = \mathbb{E} \left\{ \left( \prod_{j=1}^p \prod_{i=1}^{a-1} \mathbb{1}_{\{x_1^i \in S_j(\Delta_i)\}} \cdots \mathbb{1}_{\{x_{k_i}^i \in S_j(\Delta_i)\}} \right) \sum_{\bar{x}_a} \prod_{j=1}^p \phi_a(\bar{x}_a - \bar{S}_j^a(n - n_a)) \right\} \\ &\quad \le \mathbb{E} \left\{ \left( \prod_{j=1}^p \prod_{i=1}^{a-1} \mathbb{1}_{\{x_1^i \in S_j(\Delta_i)\}} \cdots \mathbb{1}_{\{x_{k_i}^i \in S_j(\Delta_i)\}} \right) \right. \\ &\qquad \quad \left. \times \prod_{j=1}^p \left( \sum_{\bar{x}_a} \phi_a^p(\bar{x}_a - \bar{S}_j^a(n - n_a)) \right)^{1/p} \right\} \\ &\quad = \mathbb{E} \left\{ \left( \prod_{j=1}^p \prod_{i=1}^{a-1} \mathbb{1}_{\{x_1^i \in S_j(\Delta_i)\}} \cdots \mathbb{1}_{\{x_{k_i}^i \in S_j(\Delta_i)\}} \right) \sum_{\bar{x}_a} \phi_a^p(\bar{x}_a) \right\} \\ &\quad = \left\{ \mathbb{E} \prod_{i=1}^{a-1} \mathbb{1}_{\{x_1^i \in S(\Delta_i)\}} \cdots \mathbb{1}_{\{x_{k_i}^i \in S(\Delta_i)\}} \right\}^p \cdot \mathbb{E} J_{n_a}^{k_a}. \end{aligned}$$

So we have

$$\sum_{\bar{x}_1} \cdots \sum_{\bar{x}_a} \left[ \mathbb{E} \prod_{i=1}^a \mathbb{1}_{\{x_1^i \in S(\Delta_i)\}} \cdots \mathbb{1}_{\{x_{k_i}^i \in S(\Delta_i)\}} \right]^p$$

$$\leq \mathbb{E}J_{n_a}^{k_a} \cdot \sum_{\bar{x}_1} \cdots \sum_{\bar{x}_{a-1}} \left[ \mathbb{E} \prod_{i=1}^{a-1} \mathbb{1}_{\{x_1^i \in S(\Delta_i)\}} \cdots \mathbb{1}_{\{x_{k_i}^i \in S(\Delta_i)\}} \right]^p.$$

Repeating this procedure gives (3.5). □

Immediately, we have:

Corollary 1. *For any integer $m \geq 1$,*

$$(\mathbb{E}J_{n_1+\cdots+n_a}^m)^{1/p} \leq \sum_{\substack{k_1+\cdots+k_m=m \\ k_1,\ldots,k_m \geq 0}} \frac{m!}{k_1! \cdots k_a!} (\mathbb{E}J_{n_1}^{k_1})^{1/p} \cdots (\mathbb{E}J_{n_a}^{k_a})^{1/p}.$$

*Consequently, for any $\lambda > 0$,*

$$\sum_{m=0}^{\infty} \frac{\lambda^m}{m!} (\mathbb{E}J_{n_1+\cdots+n_a}^m)^{1/p} \leq \prod_{i=1}^{a} \sum_{m=0}^{\infty} \frac{\lambda^m}{m!} (\mathbb{E}J_{n_i}^m)^{1/p}.$$

As application, we have the following sharp moment estimate.

Lemma 1. *There is a constant $C > 0$ depending only on $d$ and $p$ such that:*

(i) *When $d = 2$ and $p \geq 2$,*

$$\mathbb{E}J_n^m \leq (m!)^{p-1} C^m n^m \left( \min\left\{ \frac{1}{(\log(n/m))^p}, 1 \right\} \right)^m \qquad \forall m, n = 1, 2, \ldots. \tag{3.6}$$

(ii) *When $d = 3$, $p = 2$,*

$$\mathbb{E}J_n^m \leq (m!)^{3/2} C^m n^{m/2} \qquad \forall m, n = 1, 2, \ldots. \tag{3.7}$$

Proof. Due to similarity we only prove (3.6) in the case $\log(n/m) \geq 1$. Write $l(m,n) = [n/m] + 1$. Then

$$\begin{aligned}
(\mathbb{E}J_n^m)^{1/p} &\leq \sum_{\substack{k_1+\cdots+k_m=m \\ k_1,\ldots,k_m \geq 0}} \frac{m!}{k_1! \cdots k_m!} (\mathbb{E}J_{l(m,n)}^{k_1})^{1/p} \cdots (\mathbb{E}J_{l(m,n)}^{k_m})^{1/p} \\
&\leq \sum_{\substack{k_1+\cdots+k_m=m \\ k_1,\ldots,k_m \geq 0}} \frac{m!}{k_1! \cdots k_m!} k_1! \cdots k_m! (\mathbb{E}J_{l(m,n)})^{k_1/p} \cdots (\mathbb{E}J_{l(m,n)})^{k_m/p} \\
&= \binom{2m-1}{m} m! C^m (\mathbb{E}J_{l(m,n)})^{m/p} \\
&= \binom{2m-1}{m} m! C^m \left( \frac{(n/m)}{(\log(n/m))^p} \right)^{m/p} \\
&\leq \binom{2m}{m} (m!)^{(p-1)/p} C^m \left( \frac{n}{(\log(n/m))^p} \right)^{m/p}
\end{aligned}$$

where the second inequality follows from the fact [Remarks, page 664 in Le Gall and Rosen (1991)] that

$$\mathbb{E}J_n^k \le (k!)^p (\mathbb{E}J_n)^k, \qquad k = 0, 1, \ldots.$$

Hence

$$\mathbb{E}J_n^m \le \binom{2m}{m}^p C^{pm} (m!)^{p-1} \left(\frac{n}{(\log(n/m))^p}\right)^m.$$

Finally, the desired conclusion follows from the fact

$$\binom{2m}{m} \le 4^m. \qquad \square$$

We are ready to prove the upper bounds for Theorems 1 and 2. Due to similarity we only prove (3.1). Let $t > 0$ be fixed and let $t_n = [tn/b_n]$. Applying Corollary 1, we have

$$\begin{aligned}
&\sum_{m=0}^\infty \frac{1}{m!}\theta^m \left(\frac{b_n \log^p n}{n}\right)^{m/p} (\mathbb{E}J_n^m)^{1/p} \\
&\quad \le \left(\sum_{m=0}^\infty \frac{1}{m!}\theta^m \left(\frac{b_n \log^p n}{n}\right)^{m/p} (\mathbb{E}J_{t_n}^m)^{1/p}\right)^{[n/t_n]+1}.
\end{aligned} \tag{3.8}$$

By (1.4), Lemma 1 and the dominated convergence theorem,

$$\begin{aligned}
&\sum_{m=0}^\infty \frac{1}{m!}\theta^m \left(\frac{b_n \log^p n}{n}\right)^{m/p} (\mathbb{E}J_{t_n}^m)^{1/p} \\
&\quad \longrightarrow \sum_{m=0}^\infty \frac{1}{m!}(2\pi\theta)^m t^{m/p} \det(\Gamma)^{1/(2p)m} (\mathbb{E}\alpha([0,1]^p)^m)^{1/p}
\end{aligned} \tag{3.9}$$

as $n \to \infty$. Hence,

$$\begin{aligned}
&\limsup_{n\to\infty} \frac{1}{b_n} \log\left(\sum_{m=0}^\infty \frac{1}{m!}\theta^m \left(\frac{b_n \log^p n}{n}\right)^{m/p} (\mathbb{E}J_n^m)^{1/p}\right) \\
&\quad \le \frac{1}{t}\log\left(\sum_{m=0}^\infty \frac{1}{m!}(2\pi\theta)^m t^{m/p} \det(\Gamma)^{1/(2p)m} (\mathbb{E}\alpha([0,1]^p)^m)^{1/p}\right).
\end{aligned} \tag{3.10}$$

In view of (1.7) (with $d = 2$), applying Theorem 5 to $\varepsilon = t^{-1}$,

$$Z_\varepsilon = t^{-(p-1)}\alpha([0,1]^p) \quad \text{and} \quad I(\lambda) = \frac{p}{2}\kappa(2,p)^{-2p/(p-1)}\lambda^{1/(p-1)}$$

gives

$$\lim_{t\to\infty}\frac{1}{t}\log\left(\sum_{m=0}^{\infty}\frac{1}{m!}(2\pi\theta)^m t^{m/p}\det(\Gamma)^{1/(2p)m}(\mathbb{E}\alpha([0,1]^p)^m)^{1/p}\right)$$
$$=\sup_{\lambda>0}\left\{(2\pi\theta)\det(\Gamma)^{1/(2p)}\lambda^{1/p}-\frac{1}{2}\kappa(2,p)^{-2p/(p-1)}\lambda^{1/(p-1)}\right\}$$
$$=\frac{1}{p}\left(\frac{2(p-1)}{p}\right)^{p-1}(2\pi\theta)^p\sqrt{\det(\Gamma)}\kappa(2,p)^{2p}.$$

Letting $t\to\infty$ in (3.10) gives (3.1).

**4. Lower bounds.** The main goal of this section is to prove that in the case $d=2$, $p\geq 2$,

$$\liminf_{n\to\infty}\frac{1}{b_n}\log\sum_{m=0}^{\infty}\frac{\theta^m}{m!}\left(\frac{b_n\log^p n}{n}\right)^{m/p}(\mathbb{E}J_n^m)^{1/p}$$
$$\geq\frac{1}{p}\left(\frac{2(p-1)}{p}\right)^{p-1}(2\pi\theta)^p\sqrt{\det(\Gamma)}\kappa(2,p)^{2p}\qquad(\theta>0) \tag{4.1}$$

for any $\{b_n\}$ satisfying (1.10); and that in the case $d=3$, $p=2$,

$$\liminf_{n\to\infty}\frac{1}{b_n}\log\sum_{m=0}^{\infty}\frac{\theta^m}{m!}\left(\frac{b_n}{n}\right)^{m/4}\sqrt{\mathbb{E}J_n^m}$$
$$\geq 2\left(\frac{3}{4}\right)^3(\gamma(S)\theta)^4\det(\Gamma)^{-1}\kappa(3,2)^8\qquad(\theta>0) \tag{4.2}$$

for any $\{b_n\}$ satisfying (1.12).

We proceed in two steps. The main result in the first step is a weak law given in Theorem 7 and the essential tool is the second moment estimate. The second step starts after the proof of Theorem 7 and the goal is to establish Theorem 8 which leads to (4.1) and (4.2) through a simple argument. To this end we first establish a Feynman–Kac lower bound in Lemma 5, using an argument similar to the one given in the proof of Theorem 4.1 of Chen and Li (2004). The accomplishment of Theorem 8 relies on eliminating the contribution from self-intersection between different time periods. This part is carried out in Lemma 6.

For any $x=(x_1,\dots,x_d)\in\mathbb{R}^d$, we adopt the notation $[x]\in\mathbb{Z}^d$ throughout this section for the lattice part of $x$, that is,

$$[x]=([x_1],\dots,[x_d]).$$

Recall that a $\mathbb{Z}^d$ random walk $\{S(n)\}$ is said to be aperiodic if the greatest common factor of the set

$$\{n\geq 1;\mathbb{P}\{S(n)=0\}>0\}$$

is 1. According to a remark made in page 661 of Le Gall and Rosen (1991), the aperiodicity implies

$$\sup_{x\in\mathbb{Z}^d}\left|n^{d/2}\mathbb{P}\{S(n)=x\}-\frac{1}{(2\pi)^{d/2}\det(\Gamma)^{1/2}}\exp\left\{-\frac{1}{2n}\langle x,\Gamma^{-1}x\rangle\right\}\right|\to 0 \tag{4.3}$$

as $n\to\infty$.

Lemma 2. *Let $\{S(n)\}$ be a mean zero, square integrable random walk in $\mathbb{Z}^d$. For any $x\in\mathbb{Z}^d$, write*

$$T_x=\inf\{n\ge 0; S(n)=x\}.$$

*Then*

$$\mathbb{P}\{T_x\le n\}\ge\sum_{k=0}^n\mathbb{P}\{S(k)=x\}\Big/\sum_{k=0}^n\mathbb{P}\{S(k)=0\},\qquad n=1,2,\dots. \tag{4.4}$$

Proof. By the Markov property,

$$\begin{aligned}\mathbb{P}\{S(k)=x\}&=\sum_{j=0}^k\mathbb{P}\{T_x=j,S(k)=x\}\\&=\sum_{j=0}^k\mathbb{P}\{T_x=j\}\mathbb{P}\{S(k-j)=0\}.\end{aligned} \tag{4.5}$$

Summing up on both sides,

$$\begin{aligned}\sum_{k=0}^n\mathbb{P}\{S(k)=x\}&=\sum_{k=0}^n\sum_{j=0}^k\mathbb{P}\{T_x=j\}\mathbb{P}\{S(k-j)=0\}\\&=\sum_{j=0}^n\mathbb{P}\{T_x=j\}\sum_{k=j}^n\mathbb{P}\{S(k-j)=0\}\\&\le\mathbb{P}\{T_x\le n\}\sum_{k=0}^n\mathbb{P}\{S(k)=0\}.\end{aligned}$$

□

Lemma 3. *Let $\{S(n)\}$ be a mean zero, square integrable random walk in $\mathbb{Z}^d$.*

(i) *As $d=2$,*

$$\sup_n\mathbb{E}\exp\left\{\theta\frac{\log n}{n}\#\{S[0,n]\}\right\}<\infty\qquad(\theta>0). \tag{4.6}$$

(ii) *As $d\ge 3$,*

$$\sup_n\mathbb{E}\exp\left\{\theta\frac{1}{n}\#\{S[0,n]\}\right\}<\infty\qquad(\theta>0). \tag{4.7}$$

PROOF. Since $\#\{S[0,n]\} \le n+1$, (4.7) is trivial. To prove (4.6), we first show that for any $a, b > 0$ and any integer $n \ge 1$,

$$\mathbb{P}\{\#\{S[0,n]\} \ge a+b\} \le \mathbb{P}\{\#\{S[0,n]\} \ge a\}\mathbb{P}\{\#\{S[0,n]\} \ge b\}. \tag{4.8}$$

Notice that $\#\{S[0,n]\}$ takes integer values. So we may assume that $a$ and $b$ are integers, for otherwise we can use $[a]$, $[b]$ and $[a+b]$ instead of $a$, $b$ and $a+b$, respectively, in the following argument. Define the stopping time

$$\tau = \inf\{k; \#\{S[0,n]\} \ge a\}.$$

Then

$$\begin{aligned}
&\mathbb{P}\{\#\{S[0,n]\} \ge a+b\} \\
&\quad = \mathbb{P}\{\#\{S[0,n]\} \ge a+b, \tau \le n\} \\
&\quad = \sum_{k=0}^{n} \mathbb{P}\{\tau = k, \#\{S[0,n]\} - \#\{S[0,k]\} \ge b\} \\
&\quad \le \sum_{k=0}^{n} \mathbb{P}\{\tau = k, \#\{S[k,n]\} \ge b\} \\
&\quad = \sum_{k=0}^{n} \mathbb{P}\{\tau = k\}\mathbb{P}\{\#\{S[0,n-k]\} \ge b\} \\
&\quad \le \mathbb{P}\{\tau \le n\}\mathbb{P}\{\#\{S[0,n]\} \ge b\} \\
&\quad = \mathbb{P}\{\#\{S[0,n]\} \ge a\}\mathbb{P}\{\#\{S[0,n]\} \ge b\}.
\end{aligned}$$

We now prove (4.6) in the case $d = 2$. Let $C > 0$ be fixed. By (4.8) we have

$$\mathbb{P}\left\{\#\{S[0,n]\} \ge Cm\frac{n}{\log n}\right\} \le \left(\mathbb{P}\left\{\#\{S[0,n]\} \ge C\frac{n}{\log n}\right\}\right)^m.$$

By the fact that $\mathbb{E}\#\{S[0,n]\} = O(n(\log n)^{-1})$ one can take $C > 0$ large enough so

$$\sup_n P\left\{\#\{S[0,n]\} \ge C\frac{n}{\log n}\right\} \le e^{-2}.$$

Therefore, (4.6) holds for $\theta = C^{-1}$. We now show that it holds for all $\theta > 0$. Indeed, take $\delta > 0$ such that $\theta < C^{-1}[\delta^{-1}]$ and write $k_n = [\delta n]$. The desired conclusion follows from the following estimate:

$$\begin{aligned}
\mathbb{E}\exp\left\{\theta\frac{\log n}{n}\#\{S[0,n]\}\right\} &\le \left(\mathbb{E}\exp\left\{\theta\frac{\log n}{n}\#\{S[0,k_n]\}\right\}\right)^{[\delta^{-1}]+1} \\
&\le \left(\mathbb{E}\exp\left\{C^{-1}\frac{\log k_n}{k_n}\#\{S[0,k_n]\}\right\}\right)^{[\delta^{-1}]+1}. \qquad \square
\end{aligned}$$

THEOREM 7. *Let $\{S(n)\}$ be a mean zero, square integrable random walk in $\mathbb{Z}^d$ and let $X_t$ be the symmetric Lévy Gaussian process such that $S(1)$ and $X_1$ have the same covariance matrix $\Gamma$. Let $f(x)$ be a bounded, continuous function on $\mathbb{R}^d$.*

(i) *As $d=2$,*

$$\left(\frac{\log n}{2\pi n\sqrt{\det(\Gamma)}}\sum_{x\in S[0,n]} f\left(\frac{x}{\sqrt{n}}\right), \frac{S(n)}{\sqrt{n}}\right) \xrightarrow{d} \left(\int_0^1 f(X_t)\,dt, X_1\right). \tag{4.9}$$

(ii) *As $d\geq 3$,*

$$\left(\frac{1}{\gamma(S)n}\sum_{x\in S[0,n]} f\left(\frac{x}{\sqrt{n}}\right), \frac{S(n)}{\sqrt{n}}\right) \xrightarrow{d} \left(\int_0^1 f(X_t)\,dt, X_1\right). \tag{4.10}$$

PROOF. We only consider the case $d=2$, as the proof for $d\geq 3$ is similar. By the invariance principle,

$$\left(\frac{1}{n}\sum_{k=1}^n f\left(\frac{S(k)}{\sqrt{n}}\right), \frac{S(n)}{\sqrt{n}}\right) \xrightarrow{d} \left(\int_0^1 f(X_t)\,dt, X_1\right).$$

Let

$$l(n,x)=\sum_{k=1}^n \mathbb{1}_{\{S(k)=x\}}, \qquad x\in\mathbb{Z}^2,\ n\geq 1,$$

be the local time of $\{S(n)\}$. By the fact

$$\sum_{k=1}^n f\left(\frac{S(k)}{\sqrt{n}}\right)=\sum_{x\in\mathbb{Z}^2} f\left(\frac{x}{\sqrt{n}}\right)l(n,x)$$

we need only to prove

$$\begin{aligned}\frac{1}{n^2}\mathbb{E}\Bigg[&\sum_{x\in\mathbb{Z}^2} f\left(\frac{x}{\sqrt{n}}\right)l(n,x)\\ &-\frac{\log n}{2\pi\sqrt{\det(\Gamma)}}\sum_{x\in\mathbb{Z}^2} f\left(\frac{x}{\sqrt{n}}\right)\mathbb{1}_{\{T_x\leq n\}}\Bigg]^2\to 0\\ &\qquad (n\to\infty)\end{aligned} \tag{4.11}$$

where $T_x=\inf\{n\geq 0, S(n)=x\}$.

We may assume that $f\geq 0$, for otherwise we consider the decomposition $f=f^+-f^-$. We only need to prove

$$\frac{1}{n^2}\mathbb{E}\left[\sum_{x\in\mathbb{Z}^2} f\left(\frac{x}{\sqrt{n}}\right)l(n,x)\right]^2 \longrightarrow \mathbb{E}\left[\int_0^1 f(X_t)\,dt\right]^2 \qquad (n\to\infty), \tag{4.12}$$

$$\frac{\log^2 n}{4\pi^2 n^2 \det(\Gamma)}\mathbb{E}\bigg[\sum_{x\in\mathbb{Z}^2} f\Big(\frac{x}{\sqrt{n}}\Big)\mathbb{1}_{\{T_x\le n\}}\bigg]^2 \longrightarrow \mathbb{E}\bigg[\int_0^1 f(X_t)\,dt\bigg]^2 \qquad (n\to\infty), \tag{4.13}$$

$$\frac{\log n}{2\pi n^2\sqrt{\det(\Gamma)}}\mathbb{E}\bigg[\bigg(\sum_{x\in\mathbb{Z}^2} f\Big(\frac{x}{\sqrt{n}}\Big)\mathbb{1}_{\{T_x\le n\}}\bigg)\bigg(\sum_{x\in\mathbb{Z}^2} f\Big(\frac{x}{\sqrt{n}}\Big)l(n,x)\bigg)\bigg] \longrightarrow \mathbb{E}\bigg[\int_0^1 f(X_t)\,dt\bigg]^2 \qquad (n\to\infty). \tag{4.14}$$

Clearly, (4.12) is a direct consequence of the invariance principle and the dominated convergence theorem. Notice that

$$\begin{aligned}\sum_{x\in\mathbb{Z}^2} f\Big(\frac{x}{\sqrt{n}}\Big)\mathbb{1}_{\{T_x\le n\}} &= \int_{\mathbb{R}^2} f\Big(\frac{[x]}{\sqrt{n}}\Big)\mathbb{1}_{\{T_{[x]}\le n\}}\,dx \\ &= n\int_{\mathbb{R}^2} f\Big(\frac{[\sqrt{n}x]}{\sqrt{n}}\Big)\mathbb{1}_{\{T_{[\sqrt{n}x]}\le n\}}\,dx \\ &= o(1)\cdot\#\{S[0,n]\} + n\int_{\mathbb{R}^2} f(x)\mathbb{1}_{\{T_{[\sqrt{n}x]}\le n\}}\,dx \qquad (n\to\infty).\end{aligned}$$

By Lemma 3, (4.13) is equivalent to

$$\frac{\log^2 n}{4\pi^2 \det(\Gamma)}\mathbb{E}\bigg[\int_{\mathbb{R}^2} f(x)\mathbb{1}_{\{T_{[\sqrt{n}x]}\le n\}}\,dx\bigg]^2 \longrightarrow \mathbb{E}\bigg[\int_0^1 f(X_t)\,dt\bigg]^2 \qquad (n\to\infty).$$

Notice that

$$\mathbb{E}\bigg[\int_{\mathbb{R}^2} f(x)\mathbb{1}_{\{T_{[\sqrt{n}x]}\le n\}}\,dx\bigg]^2 = 2\int_{\mathbb{R}^2\times\mathbb{R}^2} f(x)f(y)\mathbb{P}\{T_{[\sqrt{n}x]}\le T_{[\sqrt{n}y]}\le n\}\,dx\,dy.$$

By (5.d) and (5.e) in Le Gall (1986a), respectively,

$$\lim_{n\to\infty}(\log n)^2\mathbb{P}\{T_{[\sqrt{n}x]}\le T_{[\sqrt{n}y]}\le n\} = (2\pi)^2\det(\Gamma)\iint_{\{0\le s\le t\le 1\}} p_s(x)p_{t-s}(y-x)\,ds\,dt,$$

$$(\log n)^2\mathbb{P}\{T_{[\sqrt{n}x]}\le T_{[\sqrt{n}y]}\le n\} \le C^2 h(|x|)h(|y-x|),$$

where $p_t(x)$ is the density of $X_t$ and $h(r) = (\log(1/r))_+ + r^{-2}\mathbb{1}_{\{r>1/2\}}$. By the dominated convergence theorem,

$$\lim_{n\to\infty} \frac{\log^2 n}{4\pi^2 \det(\Gamma)} \mathbb{E}\Big[\int_{\mathbb{R}^2} f(x)\mathbb{1}_{\{T_{[\sqrt{n}x]}\le n\}}\,dx\Big]^2$$

$$= 2\int_{\mathbb{R}^2\times\mathbb{R}^2} f(x)f(y)\Big\{\iint_{\{0\le s\le t\le 1\}} p_s(x)p_{t-s}(y-x)\,ds\,dt\Big\}\,dx\,dy$$

$$= 2\iint_{\{0\le s\le t\le 1\}} ds\,dt \int_{\mathbb{R}^2} dx\, f(x)p_s(x)\int_{\mathbb{R}^2} f(y)p_{t-s}(y-x)\,dy$$

$$= 2\iint_{\{0\le s\le t\le 1\}} ds\,dt\,\mathbb{E}\{f(X_s)\mathbb{E}_{X_s} f(X_{t-s})\}$$

$$= 2\iint_{\{0\le s\le t\le 1\}} \mathbb{E}\{f(X_s)f(X_t)\}\,ds\,dt = \mathbb{E}\Big[\int_0^1 f(X_t)\,dt\Big]^2.$$

We now come to the proof of (4.14). Since

$$\mathbb{E}\Big[\Big(\sum_{x\in\mathbb{Z}^2} f\Big(\frac{x}{\sqrt{n}}\Big)\mathbb{1}_{\{T_x\le n\}}\Big)\Big(\sum_{x\in\mathbb{Z}^2} f\Big(\frac{x}{\sqrt{n}}\Big)l(n,x)\Big)\Big]$$

$$\le \Big\{\mathbb{E}\Big[\sum_{x\in\mathbb{Z}^2} f\Big(\frac{x}{\sqrt{n}}\Big)\mathbb{1}_{\{T_x\le n\}}\Big]^2\Big\}^{1/2}\Big\{\mathbb{E}\Big[\sum_{x\in\mathbb{Z}^2} f\Big(\frac{x}{\sqrt{n}}\Big)l(n,x)\Big]^2\Big\}^{1/2},$$

by (4.12) and (4.13)

$$\limsup_{n\to\infty} \frac{\log n}{2\pi n^2\sqrt{\det(\Gamma)}}\mathbb{E}\Big[\Big(\sum_{x\in\mathbb{Z}^2} f\Big(\frac{x}{\sqrt{n}}\Big)\mathbb{1}_{\{T_x\le n\}}\Big)\Big(\sum_{x\in\mathbb{Z}^2} f\Big(\frac{x}{\sqrt{n}}\Big)l(n,x)\Big)\Big]$$

$$\le \mathbb{E}\Big[\int_0^1 f(X_t)\,dt\Big]^2.$$

To obtain the lower bound for (4.14), notice that

$$\mathbb{E}\Big[\Big(\sum_{x\in\mathbb{Z}^2} f\Big(\frac{x}{\sqrt{n}}\Big)\mathbb{1}_{\{T_x\le n\}}\Big)\Big(\sum_{x\in\mathbb{Z}^2} f\Big(\frac{x}{\sqrt{n}}\Big)l(n,x)\Big)\Big]$$

$$\sim \sum_{x,y\in\mathbb{Z}^2} f\Big(\frac{x}{\sqrt{n}}\Big)f\Big(\frac{y}{\sqrt{n}}\Big)\sum_{0\le j\le k\le n}\mathbb{P}\{T_x=j, S(k)=y\}$$

$$+ \sum_{x,y\in\mathbb{Z}^2} f\Big(\frac{x}{\sqrt{n}}\Big)f\Big(\frac{y}{\sqrt{n}}\Big)\sum_{0\le j\le k\le n}\mathbb{P}\{S(j)=x, T_y=k\}.$$

By the Markov property,

$$\begin{aligned}
&\sum_{0\le j\le k\le n} \mathbb{P}\{T_x = j, S(k) = y\} \\
&\quad = \sum_{0\le j\le k\le n} \mathbb{P}\{T_x = j\}\mathbb{P}\{S(k-j) = y-x\} \\
&\quad = \sum_{k=1}^{n} \mathbb{P}\{S(k) = y-x\}\mathbb{P}\{T_x \le n-k\} \\
&\quad \ge \sum_{k=1}^{n} \mathbb{P}\{S(k) = y-x\}(G(n-k))^{-1}\sum_{j=1}^{n-k} \mathbb{P}\{S(j) = x\} \\
&\quad \ge (G(n))^{-1}\sum_{k=1}^{n} \mathbb{P}\{S(k) = y-x\}\sum_{j=1}^{n-k} \mathbb{P}\{S(j) = x\} \\
&\quad = (G(n))^{-1}\sum_{0\le j\le k\le n} \mathbb{P}\{S(j) = x\}\mathbb{P}\{S(k-j) = y-x\},
\end{aligned}$$

where, by Proposition 2.4 in Le Gall and Rosen (1991),

$$G(n) \equiv \sum_{k=0}^{n} \mathbb{P}\{S(k) = 0\} \sim \frac{1}{2\pi\sqrt{\det(\Gamma)}}\log n \qquad (n \to \infty)$$

and where the third step follows from Lemma 2.

Using the Markov property again,

$$\begin{aligned}
&\sum_{0\le j\le k\le n} \mathbb{P}\{S(j) = x, T_y = k\} \\
&\quad = \sum_{0\le j\le k\le n} \mathbb{P}\{S(j) = x, T_y \ge j, S(j) \ne y, \dots, S(k-1) \ne y, S(k) = y\} \\
&\quad = \sum_{0\le j\le k\le n} \mathbb{P}\{S(j) = x, T_y \ge j\}\mathbb{P}\{T_{y-x} = k-j\} \\
&\quad = \sum_{0\le j\le k\le n} \mathbb{P}\{S(j) = x\}\mathbb{P}\{T_{y-x} = k-j\} \\
&\qquad - \sum_{0\le j\le k\le n} \mathbb{P}\{S(j) = x, T_y < j\}\mathbb{P}\{T_{y-x} = k-j\}.
\end{aligned}$$

For the first term on the right-hand side,

$$\begin{aligned}
&\sum_{0\le j\le k\le n} \mathbb{P}\{S(j)=x\}\mathbb{P}\{T_{y-x}=k-j\} \\
&\quad = \sum_{j=0}^{n} \mathbb{P}\{S(j)=x\}\mathbb{P}\{T_{y-x}\le n-j\} \\
&\quad \ge \sum_{j=0}^{n} \mathbb{P}\{S(j)=x\}(G(n-j))^{-1}\sum_{k=0}^{n-j}\mathbb{P}\{S(k)=y-x\} \\
&\quad \ge (G(n))^{-1}\sum_{0\le j\le k\le n}\mathbb{P}\{S(j)=x\}\mathbb{P}\{S(k-j)=y-x\}.
\end{aligned}$$

For the second term,

$$\begin{aligned}
&\sum_{0\le j\le k\le n} \mathbb{P}\{S(j)=x, T_y<j\}\mathbb{P}\{T_{y-x}=k-j\} \\
&\quad \le \mathbb{P}\{T_{y-x}\le n\}\sum_{j=0}^{n}\mathbb{P}\{S(j)=x,T_y<j\} \\
&\quad = \mathbb{P}\{T_{y-x}\le n\}\sum_{j=0}^{n}\sum_{i=0}^{j}\mathbb{P}\{T_y=i,S(j)=x\} \\
&\quad = \mathbb{P}\{T_{y-x}\le n\}\sum_{j=0}^{n}\sum_{i=0}^{j}\mathbb{P}\{T_y=i\}\mathbb{P}\{S(j-i)=x-y\} \\
&\quad \le \mathbb{P}\{T_x\le n\}\mathbb{P}\{T_{y-x}\le n\}\sum_{j=0}^{n}\mathbb{P}\{S(j)=x-y\}.
\end{aligned}$$

Summarizing what we have,

$$\begin{aligned}
&\liminf_{n\to\infty}\frac{\log n}{2\pi n^2\sqrt{\det(\Gamma)}}\mathbb{E}\left[\left(\sum_{x\in\mathbb{Z}^2} f\left(\frac{x}{\sqrt{n}}\right)\mathbb{1}_{\{T_x\le n\}}\right)\left(\sum_{x\in\mathbb{Z}^2} f\left(\frac{x}{\sqrt{n}}\right)l(n,x)\right)\right] \\
&\quad \ge \liminf_{n\to\infty}\frac{1}{n^2}\mathbb{E}\left[\sum_{x\in\mathbb{Z}^2} f\left(\frac{x}{\sqrt{n}}\right)l(n,x)\right]^2 \\
&\qquad - \limsup_{n\to\infty}\frac{\log n}{2\pi n^2\sqrt{\det(\Gamma)}}\sum_{x,y\in\mathbb{Z}^2} f\left(\frac{x}{\sqrt{n}}\right)f\left(\frac{y}{\sqrt{n}}\right) \\
&\qquad\qquad \times \mathbb{P}\{T_x\le n\}\mathbb{P}\{T_{y-x}\le n\}\sum_{j=0}^{n}\mathbb{P}\{S(j)=x-y\}.
\end{aligned}$$

In view of (4.12), it remains to prove

$$\lim_{n\to\infty}\frac{\log n}{n^2}\sum_{x,y\in\mathbb{Z}^2} f\Big(\frac{x}{\sqrt{n}}\Big)f\Big(\frac{y}{\sqrt{n}}\Big) \times\mathbb{P}\{T_x\le n\}\mathbb{P}\{T_{y-x}\le n\}\sum_{j=0}^n\mathbb{P}\{S(j)=x-y\}=0. \tag{4.15}$$

Indeed,

$$\begin{aligned}&\sum_{x,y\in\mathbb{Z}^2} f\Big(\frac{x}{\sqrt{n}}\Big)f\Big(\frac{y}{\sqrt{n}}\Big)\mathbb{P}\{T_x\le n\}\mathbb{P}\{T_{y-x}\le n\}\sum_{j=0}^n\mathbb{P}\{S(j)=x-y\}\\ &\quad\le \|f\|_\infty\sum_{x,y\in\mathbb{Z}^2} f\Big(\frac{x}{\sqrt{n}}\Big)\mathbb{P}\{T_x\le n\}\mathbb{P}\{T_{y-x}\le n\}\sum_{j=0}^n\mathbb{P}\{S(j)=x-y\}\\ &\quad\le \|f\|_\infty\bigg\{\mathbb{E}\sum_{x\in\mathbb{Z}^2} f\Big(\frac{x}{\sqrt{n}}\Big)\mathbb{1}_{\{T_x\le n\}}\bigg\}\bigg\{\sum_{x\in\mathbb{Z}^2}\mathbb{P}\{T_x\le n\}\sum_{j=0}^n\mathbb{P}\{S(j)=-x\}\bigg\}.\end{aligned}$$

From (4.13),

$$\limsup_{n\to\infty}\frac{\log n}{n}\mathbb{E}\sum_{x\in\mathbb{Z}^2} f\Big(\frac{x}{\sqrt{n}}\Big)\mathbb{1}_{\{T_x\le n\}}<\infty.$$

Notice that

$$\begin{aligned}&\sum_{x\in\mathbb{Z}^2}\mathbb{P}\{T_x\le n\}\sum_{j=0}^n\mathbb{P}\{S(j)=-x\}\\ &\quad\le\bigg\{\sum_{x\in\mathbb{Z}^2}(\mathbb{P}\{T_x\le n\})^2\bigg\}^{1/2}\bigg\{\sum_{x\in\mathbb{Z}^2}\bigg[\sum_{j=0}^n\mathbb{P}\{S(j)=x\}\bigg]^2\bigg\}^{1/2}.\end{aligned}$$

Finally, (4.15) follows from the fact that as $p=2$,

$$\mathbb{E}J_n=\sum_{x\in\mathbb{Z}^2}(\mathbb{P}\{T_x\le n\})^2=O\Big(\frac{n}{(\log n)^2}\Big),$$

$$\mathbb{E}I_n=\sum_{x\in\mathbb{Z}^2}\bigg[\sum_{j=0}^n\mathbb{P}\{S(j)=x\}\bigg]^2=O(n).$$

□

Fix integer $t\ge 1$ and the bounded measurable function $f$ on $\mathbb{R}^d$. Define the linear operator $T$ on $\mathcal{L}^2(\mathbb{Z}^d)$ by

$$\begin{aligned}(T\xi)(x)&=\mathbb{E}_x\bigg[\exp\bigg\{\sum_{y\in S[0,t]}f(y)\bigg\}\xi(S(t))\bigg]\\ &=\mathbb{E}\bigg[\exp\bigg\{\sum_{y\in S[0,t]}f(x+y)\bigg\}\xi(x+S(t))\bigg].\end{aligned}$$

LEMMA 4. *Given any symmetric random walk* $\{S(n)\}$ *on* $\mathbb{Z}^d$, $T$ *is self-adjoint: For any* $\xi,\eta\in\mathcal{L}^2(\mathbb{Z}^d)$, $\langle\eta,T\xi\rangle=\langle T\eta,\xi\rangle$.

PROOF.

$$\begin{aligned}\langle\eta,T\xi\rangle&=\sum_{x\in\mathbb{Z}^d}\eta(x)\mathbb{E}\left[\exp\left\{\sum_{y\in S[0,t]}f(x+y)\right\}\xi(x+S(t))\right]\\&=\mathbb{E}\left[\sum_{x\in\mathbb{Z}^d}\eta(x-S(t))\exp\left\{\sum_{y\in S[0,t]}f(x+y-S(t))\right\}\xi(x)\right]\\&=\mathbb{E}\left[\sum_{x\in\mathbb{Z}^d}\eta(x+S'(t))\exp\left\{\sum_{y\in S'[0,t]}f(x+y)\right\}\xi(x)\right]\\&=\mathbb{E}\left[\sum_{x\in\mathbb{Z}^d}\eta(x+S(t))\exp\left\{\sum_{y\in S[0,t]}f(x+y)\right\}\xi(x)\right]\\&=\langle T\eta,\xi\rangle,\end{aligned}$$

where $S'(k)=-S(t)+S(t-k)$, $k=0,1,\dots,t$ and the fourth equality follows from the fact that

$$\{S'(0),\dots,S'(t)\}\stackrel{d}{=}\{S(0),\dots,S(t)\}.\qquad\square$$

In the rest of the paper, we adopt the notation

$$t_n=[n/b_n]\quad\text{and}\quad\Delta_i=[(i-1)t_n,it_n],\qquad i=1,2,\dots.\tag{4.16}$$

Write

$$\mathcal{F}_d=\{g\in\mathcal{L}^2(\mathbb{R}^d);\ \|g\|_2=1\text{ and }\|\nabla g\|_2<\infty\}.$$

LEMMA 5. *Let* $\{S(n)\}$ *be a symmetric, square integrable and aperiodic random walk on* $\mathbb{Z}^d$ *and let* $f$ *be bounded and continuous on* $\mathbb{R}^d$. *Assume that* $\{b_n\}$ *satisfies* (1.16).

(i) *As* $d=2$,

$$\begin{aligned}&\liminf_{n\to\infty}\frac{1}{b_n}\log\mathbb{E}\exp\left\{\frac{b_n\log(n/b_n)}{2\pi n\sqrt{\det(\Gamma)}}\sum_{i=1}^{[b_n]}\sum_{x\in S(\Delta_i)}f\left(\sqrt{\frac{b_n}{n}}x\right)\right\}\\&\quad\ge\sup_{g\in\mathcal{F}_2}\left\{\int_{\mathbb{R}^2}f(x)g^2(x)\,dx-\frac12\int_{\mathbb{R}^2}\langle\nabla g(x),\Gamma\nabla g(x)\rangle\,dx\right\}.\end{aligned}\tag{4.17}$$

(ii) *As* $d \geq 3$,

$$\liminf_{n\to\infty} \frac{1}{b_n} \log \mathbb{E} \exp\left\{ \frac{b_n}{\gamma(S)n} \sum_{i=1}^{[b_n]} \sum_{x\in S(\Delta_i)} f\left(\sqrt{\frac{b_n}{n}}x\right)\right\}$$
$$\geq \sup_{g\in\mathcal{F}_3} \left\{ \int_{\mathbb{R}^3} f(x)g^2(x)\,dx - \frac{1}{2}\int_{\mathbb{R}^3} \langle \nabla g(x), \Gamma\nabla g(x)\rangle\,dx\right\}. \tag{4.18}$$

PROOF. We only consider the case $d = 2$, as the proof for $d \geq 3$ is similar. For each $n$, define the continuous, self-adjoint linear operator $T_n$ on $\mathcal{L}^2(\mathbb{Z}^2)$ as

$$T_n\xi(x) = \mathbb{E}_x\left( \exp\left\{ \frac{b_n \log(n/b_n)}{2\pi n\sqrt{\det(\Gamma)}} \sum_{x\in S[0,t_n]} f\left(\sqrt{\frac{b_n}{n}}x\right)\right\}\xi(S_{t_n})\right)$$

where $x \in \mathbb{Z}^2$ and $\xi \in \mathcal{L}^2(\mathbb{Z}^2)$.

Let $g$ be a bounded function on $\mathbb{R}^2$ and assume that $g$ is infinitely differentiable, supported by a finite box $[-M, M]^2$ and

$$\int_{\mathbb{R}^2} |g(x)|^2\,dx = 1$$

and write

$$\xi_n(x) = g\left(\sqrt{\frac{b_n}{n}}x\right)\Big/\sqrt{\sum_{y\in\mathbb{Z}^2} g^2\left(\sqrt{\frac{b_n}{n}}y\right)}, \qquad x \in \mathbb{Z}^2.$$

Let $P_{t_n}(x)$ ($x \in \mathbb{Z}^2$) be the probability density of $S_{t_n}$. Then

$$\begin{aligned}
&\mathbb{E}\exp\left\{ \frac{b_n \log(n/b_n)}{2\pi n\sqrt{\det(\Gamma)}} \sum_{i=2}^{[b_n]} \sum_{x\in S(\Delta_i)} f\left(\sqrt{\frac{b_n}{n}}x\right)\right\}\\
&\quad= \sum_{x\in\mathbb{Z}^2} P_{t_n}(x)\mathbb{E}_x \exp\left\{ \frac{b_n \log(n/b_n)}{2\pi n\sqrt{\det(\Gamma)}} \sum_{i=1}^{[b_n]-1} \sum_{x\in S(\Delta_i)} f\left(\sqrt{\frac{b_n}{n}}x\right)\right\}\\
&\quad\geq \frac{1}{\sup_y |g(y)|^2}\left\{ \sum_{y\in\mathbb{Z}^2} g^2\left(\sqrt{\frac{b_n}{n}}y\right)\right\} \cdot \sum_{x\in\mathbb{Z}^2} P_{t_n}(x)\xi_n(x)\\
&\qquad\times \mathbb{E}_x\left( \exp\left\{ \frac{b_n \log(n/b_n)}{2\pi n\sqrt{\det(\Gamma)}} \sum_{i=1}^{[b_n]-1} \sum_{x\in S(\Delta_i)} f\left(\sqrt{\frac{b_n}{n}}x\right)\right\}\xi_n(S(([b_n]-1)t_n))\right)\\
&\quad= \frac{1}{\sup_y |g(y)|^2}\left\{ \sum_{y\in\mathbb{Z}^2} g^2\left(\sqrt{\frac{b_n}{n}}y\right)\right\} \cdot \sum_{x\in\mathbb{Z}^2} P_{t_n}(x)\xi_n(x)T_n^{[b_n]-1}\xi_n(x)
\end{aligned}$$

where the last step follows from the Markov property. Notice that

$$\sum_{y\in\mathbb{Z}^2} g^2\left(\sqrt{\frac{b_n}{n}}y\right) \sim \frac{n}{b_n}\int_{\mathbb{R}^2} |g(x)|^2\,dx = \frac{n}{b_n}$$

as $n\to\infty$. In view of (4.3), by aperiodicity

$$\sup_{x\in\mathbb{Z}^2}\left|t_nP_{t_n}(x) - \frac{1}{(2\pi)\det(\Gamma)^{1/2}}\exp\left\{-\frac{1}{2t_n}\langle x,\Gamma^{-1}x\rangle\right\}\right| \to 0 \qquad (n\to\infty).$$

Since $\xi_n(x)=0$ outside $[-M\sqrt{nb_n^{-1}}, M\sqrt{nb_n^{-1}}]^2$, there is a $\delta>0$ independent of $n$, such that

$$\begin{aligned}\mathbb{E}\exp\left\{\frac{b_n\log(n/b_n)}{2\pi n\sqrt{\det(\Gamma)}}\sum_{i=2}^{[b_n]}\sum_{x\in S(\Delta_i)} f\left(\sqrt{\frac{b_n}{n}}x\right)\right\}\\ \ge \delta\sum_{x\in\mathbb{Z}^2}\xi_n(x)T_n^{[b_n]-1}\xi_n(x) = \delta\langle\xi_n, T_n^{[b_n]-1}\xi_n\rangle.\end{aligned}$$

Consider the spectral representation of $T_n$:

$$\langle\xi_n, T_n\xi_n\rangle = \int_0^\infty \lambda\mu_{\xi_n}(d\lambda)$$

where $\mu_{\xi_n}$ is a probability measure on $\mathbb{R}^+$. By the mapping theorem,

$$\begin{aligned}\langle\xi_n, T_n^{[b_n]-1}\xi_n\rangle &= \int_0^\infty \lambda^{[b_n]-1}\mu_{\xi_n}(d\lambda)\\ &\ge \left(\int_0^\infty \lambda\mu_{\xi_n}(d\lambda)\right)^{[b_n]-1} = \langle\xi_n, T_n\xi_n\rangle^{[b_n]-1}\end{aligned}$$

where the second step follows from the Jensen inequality. Hence,

$$\begin{aligned}\liminf_{n\to\infty}\frac{1}{b_n}\log\mathbb{E}\exp\left\{\frac{b_n\log(n/b_n)}{2\pi n\sqrt{\det(\Gamma)}}\sum_{i=1}^{[b_n]}\sum_{x\in S(\Delta_i)} f\left(\sqrt{\frac{b_n}{n}}x\right)\right\}\\ \ge \liminf_{n\to\infty}\log\langle\xi_n, T_n\xi_n\rangle.\end{aligned}$$

Let the Lévy Gaussian process $X_t$ be given in Theorem 7. Then

$$\begin{aligned}\langle\xi_n, T_n\xi_n\rangle = \left(\sum_{y\in\mathbb{Z}^2} g^2\left(\sqrt{\frac{b_n}{n}}y\right)\right)^{-1}\cdot\sum_{x\in\mathbb{Z}^2} g\left(\sqrt{\frac{b_n}{n}}x\right)\\ \times\mathbb{E}_x\left(\exp\left\{\frac{b_n\log(n/b_n)}{2\pi n\sqrt{\det(\Gamma)}}\sum_{y\in S[0,t_n]} f\left(\sqrt{\frac{b_n}{n}}y\right)\right\}g\left(\sqrt{\frac{b_n}{n}}S(t_n)\right)\right)\end{aligned}$$

$$
\begin{aligned}
&= (1+o(1))\left(\frac{b_n}{n}\right)\sum_{x\in\mathbb{Z}^2} g\left(\sqrt{\frac{b_n}{n}}x\right)\\
&\quad\times \mathbb{E}\left(\exp\left\{\frac{b_n\log(n/b_n)}{2\pi n\sqrt{\det(\Gamma)}}\sum_{y\in S[0,t_n]} f\left(\sqrt{\frac{b_n}{n}}(x+y)\right)\right\}\right.\\
&\qquad\qquad\qquad\left.\times g\left(\sqrt{\frac{b_n}{n}}(x+S(t_n))\right)\right)\\
&\longrightarrow \int_{\mathbb{R}^2} g(x)\mathbb{E}_x\left[\exp\left\{\int_0^1 f(X_s)\,ds\right\}g(X_1)\right]dx \qquad (n\to\infty)
\end{aligned}
$$

where the last step follows from Theorem 7, Lemma 3 and the dominated convergence theorem.

Summarizing what we have so far, we obtain

$$
\begin{aligned}
&\liminf_{n\to\infty}\frac{1}{b_n}\log\mathbb{E}\exp\left\{\frac{b_n\log(n/b_n)}{2\pi n\sqrt{\det(\Gamma)}}\sum_{i=1}^{[b_n]}\sum_{x\in S(\Delta_i)} f\left(\sqrt{\frac{b_n}{n}}x\right)\right\}\\
&\qquad\geq \log\int_{\mathbb{R}^2} g(x)\mathbb{E}_x\left[\exp\left\{\int_0^1 f(X_s)\,ds\right\}g(X_1)\right]dx.
\end{aligned}
\tag{4.19}
$$

What follows next is a standard treatment [see, e.g., Remillard (2000)] which is briefly described here: Let the semigroup of linear operators $\{\Pi_t\}$ on $\mathcal{L}^2(\mathbb{R}^2)$ be defined as

$$
\Pi_t h(x) = \mathbb{E}_x\left[\exp\left\{\int_0^t f(X_s)\,ds\right\}h(X_t)\right], \qquad h\in\mathcal{L}^2(\mathbb{R}^2),\ t\geq 0.
$$

The infinitesimal generator of $\{\Pi_t\}$ is

$$
\mathcal{A}h(x) = \frac{1}{2}\sum_{i,j=1}^{d} a_{ij}\frac{\partial^2 h}{\partial x_i\,\partial x_j}(x) + f(x)h(x)
$$

where $a_{ij}$ $(1\leq i,j\leq d)$ are entries of the matrix $\Gamma$. Clearly, $\mathcal{A}$ is self-adjoint. Let

$$
(g,\mathcal{A}g) = \int_{-\infty}^{\infty}\lambda\mu_g(d\lambda)
$$

be the spectral representation of the quadratic form $(g,\mathcal{A}g)$, where $\mu_g$ is a probability measure on $(-\infty,\infty)$. By the Jensen inequality,

$$
\int_{\mathbb{R}^2} g(x)\mathbb{E}_x\left[\exp\left\{\int_0^1 f(X_s)\,ds\right\}g(X_1)\right]dx
$$

$$\begin{aligned}
&= (g, \Pi_1 g)\\
&= \int_{-\infty}^{\infty} e^{\lambda}\mu_g(d\lambda) \geq \exp\left\{\int_{-\infty}^{\infty}\lambda\mu_g(d\lambda)\right\}\\
&= \exp\{\langle g, \mathcal{A}g\rangle\}\\
&= \exp\left\{\int_{\mathbb{R}^2} f(x)g^2(x)\,dx - \frac{1}{2}\int_{\mathbb{R}^2}\langle\nabla g(x), \Gamma\nabla g(x)\rangle\,dx\right\}.
\end{aligned}$$

In view of (4.19), taking the supremum over $g$ ends the proof. □

Recall that $t_n$ and $\Delta_i$ are defined by (4.16).

LEMMA 6. *Let $\{S(n)\}$ be a mean zero and square integrable random walk on $\mathbb{Z}^d$ and let $\varepsilon > 0$ be fixed but arbitrary.*

(i) *As $d = 2$ and $\{b_n\}$ satisfies (1.10),*

$$\limsup_{n\to\infty}\frac{1}{b_n}\log\mathbb{P}\left\{\sum_{1\leq j<k\leq[b_n]}\#\{S(\Delta_j)\cap S(\Delta_k)\} \geq \varepsilon\frac{n}{\log n}\right\} = -\infty. \tag{4.20}$$

(ii) *As $d = 3$ and $\{b_n\}$ satisfies (1.12),*

$$\limsup_{n\to\infty}\frac{1}{b_n}\log\mathbb{P}\left\{\sum_{1\leq j<k\leq[b_n]}\#\{S(\Delta_j)\cap S(\Delta_k)\} \geq \varepsilon n\right\} = -\infty. \tag{4.21}$$

PROOF. Due to similarity we only prove (4.20). To be consistent with the notation used in this paper, $\{S_1(n)\}$ and $\{S_2(n)\}$ are two independent copies of $\{S(n)\}$ and $J_n = \#\{S_1[0,n]\cap S_2[0,n]\}$. Notice that

$$\sum_{1\leq j<k\leq[b_n]}\#\{S(\Delta_j)\cap S(\Delta_k)\} = \sum_{j=1}^{[b_n]-1}\sum_{i=j+1}^{[b_n]}\#\{S(\Delta_j)\cap S(\Delta_i)\}$$

and that for any fixed $1 \leq j \leq [b_n]-1$,

$$\begin{aligned}
\sum_{i=j+1}^{[b_n]}\#\{S(\Delta_j)\cap S(\Delta_i)\} &\overset{d}{=} \sum_{i=1}^{[b_n]-j}\#\{(-S_1(\Delta_1))\cap S_2(\Delta_i)\}\\
&\leq \sum_{i=1}^{[b_n]}\#\{(-S_1(\Delta_1))\cap S_2(\Delta_i)\}.
\end{aligned}$$

By the triangular inequality, we need only to prove

$$\limsup_{n\to\infty}\frac{1}{b_n}\log\mathbb{P}\left\{\sum_{i=1}^{[b_n]}\#\{(-S_1(\Delta_1))\cap S_2(\Delta_i)\} \geq \varepsilon\frac{n}{b_n\log n}\right\} = -\infty. \tag{4.22}$$

Indeed,

$$\sum_{i=1}^{[b_n]} \#\{(-S_1(\Delta_1)) \cap S_2(\Delta_i)\} = \sum_{x\in\mathbb{Z}^2} \mathbb{1}_{\{-x\in S_1(\Delta_1)\}} \sum_{i=1}^{[b_n]} \mathbb{1}_{\{x\in S_2(\Delta_i)\}}.$$

So for any integer $m \geq 1$,

$$\begin{aligned}
&\mathbb{E}\Bigg[\sum_{i=1}^{[b_n]} \#\{(-S_1(\Delta_1)) \cap S_2(\Delta_i)\}\Bigg]^m \\
&\quad = \sum_{x_1,\dots,x_m} \Bigg[\mathbb{E}\prod_{k=1}^m \mathbb{1}_{\{-x_k\in S(\Delta_1)\}}\Bigg]\Bigg[\mathbb{E}\prod_{k=1}^m \sum_{i=1}^{[b_n]} \mathbb{1}_{\{x_k\in S(\Delta_i)\}}\Bigg] \\
&\quad \leq \Bigg\{\sum_{x_1,\dots,x_m} \Bigg[\mathbb{E}\prod_{k=1}^m \mathbb{1}_{\{x_k\in S(\Delta_1)\}}\Bigg]^2\Bigg\}^{1/2} \Bigg\{\sum_{x_1,\dots,x_m} \Bigg[\mathbb{E}\prod_{k=1}^m \sum_{i=1}^{[b_n]} \mathbb{1}_{\{x_k\in S(\Delta_i)\}}\Bigg]^2\Bigg\}^{1/2} \\
&\quad = (\mathbb{E}J_{t_n}^m)^{1/2} \Bigg\{\mathbb{E}\Bigg[\sum_{x\in\mathbb{Z}^2} \prod_{j=1}^2 \sum_{i=1}^{[b_n]} \mathbb{1}_{\{x\in S_j(\Delta_i)\}}\Bigg]^m\Bigg\}^{1/2}.
\end{aligned}$$

Hence, for any $\theta > 0$,

$$\begin{aligned}
&\sum_{m=0}^\infty \frac{\theta^m}{m!}\bigg(\frac{b_n^{3/2}(\log n)^2}{n}\bigg)^{m/2} \Bigg\{\mathbb{E}\Bigg[\sum_{i=1}^{[b_n]} \#\{(-S_1(\Delta_1)) \cap S_2(\Delta_i)\}\Bigg]^m\Bigg\}^{1/2} \\
&\quad \leq \Bigg[\sum_{m=0}^\infty \frac{\theta^m}{m!}\bigg(\frac{b_n^2(\log n)^2}{n}\bigg)^{m/2} (\mathbb{E}J_{t_n}^m)^{1/2}\Bigg]^{1/2} \\
&\qquad \times \Bigg(\sum_{m=0}^\infty \frac{\theta^m}{m!}\bigg(\frac{b_n(\log n)^2}{n}\bigg)^{m/2} \Bigg\{\mathbb{E}\Bigg[\sum_{x\in\mathbb{Z}^2} \prod_{j=1}^2 \sum_{i=1}^{[b_n]} \mathbb{1}_{\{x\in S_j(\Delta_i)\}}\Bigg]^m\Bigg\}^{1/2}\Bigg)^{1/2}.
\end{aligned}$$

Applying (3.1) with $p = 2$ and with $n$ being replaced by $t_n$, we have

$$\limsup_{n\to\infty} \frac{1}{b_n} \log \Bigg[\sum_{m=0}^\infty \frac{\theta^m}{m!}\bigg(\frac{b_n^2(\log n)^2}{n}\bigg)^{m/2} (\mathbb{E}J_{t_n}^m)^{1/2}\Bigg] \leq C_1\theta^2.$$

By (3.4) with $p = 2$,

$$\begin{aligned}
&\sum_{m=0}^\infty \frac{1}{m!}\bigg(\frac{b_n(\log n)^2}{n}\bigg)^{m/2} \Bigg\{\mathbb{E}\Bigg[\sum_{x\in\mathbb{Z}^2} \prod_{j=1}^2 \sum_{i=1}^{[b_n]} \mathbb{1}_{\{x\in S_j(\Delta_i)\}}\Bigg]^m\Bigg\}^{1/2} \\
&\qquad \leq \Bigg[\sum_{m=0}^\infty \frac{1}{m!}\bigg(\frac{b_n(\log n)^2}{n}\bigg)^{m/2} \{\mathbb{E}J_{t_n}^m\}^{1/2}\Bigg]^{b_n}.
\end{aligned}$$

Using (3.9) there is $C_2>0$ such that

$$\limsup_{n\to\infty}\frac{1}{b_n}\log\sum_{m=0}^\infty\frac{1}{m!}\Big(\frac{b_n(\log n)^2}{n}\Big)^{m/2}\times\bigg\{\mathbb{E}\bigg[\sum_{x\in\mathbb{Z}^2}\prod_{j=1}^2\sum_{i=1}^{[b_n]}\mathbb{1}_{\{x\in S_j(\Delta_i)\}}\bigg]^m\bigg\}^{1/2}\le C_2.$$

Replacing $b_n$ by $\theta^2 b_n$ gives

$$\limsup_{n\to\infty}\frac{1}{b_n}\log\sum_{m=0}^\infty\frac{\theta^m}{m!}\Big(\frac{b_n(\log n)^2}{n}\Big)^{m/2}\times\bigg\{\mathbb{E}\bigg[\sum_{x\in\mathbb{Z}^2}\prod_{j=1}^2\sum_{i=1}^{[b_n]}\mathbb{1}_{\{x\in S_j(\Delta_i)\}}\bigg]^m\bigg\}^{1/2}\le C_2\theta^2. \tag{4.23}$$

Combining the above observations there is $C_3>0$ such that for any $\theta>0$,

$$\limsup_{n\to\infty}\frac{1}{b_n}\log\sum_{m=0}^\infty\frac{\theta^m}{m!}\Big(\frac{b_n^{3/2}(\log n)^2}{n}\Big)^{m/2}\times\bigg\{\mathbb{E}\bigg[\sum_{i=1}^{[b_n]}\#\{(-S_1(\Delta_1))\cap S_2(\Delta_i)\}\bigg]^m\bigg\}^{1/2}\le C_3\theta^2.$$

Applying (2.3) in Theorem 4 we can find $\delta>0$ such that

$$\limsup_{n\to\infty}\frac{1}{b_n}\log\mathbb{P}\bigg\{\sum_{i=1}^{[b_n]}\#\{(-S_1(\Delta_1))\cap S_2(\Delta_i)\}\ge\lambda\frac{nb_n^{1/2}}{(\log n)^2}\bigg\}\le-\delta\lambda.$$

Therefore, (4.22) follows from (1.10). □

Let $p\ge 2$ be the integer given in Theorem 1 and let $q>1$ be the conjugate of $p$ defined by the relation $p^{-1}+q^{-1}=1$.

THEOREM 8. *Let $\{S(n)\}$ be a symmetric, square integrable random walk on $\mathbb{Z}^d$. Let $f$ be a nonnegative, bounded and uniformly continuous function on $\mathbb{R}^d$.*

(i) *As $d=2$, $f\in\mathcal{L}^q(\mathbb{R}^2)$ and $\{b_n\}$ satisfies (1.10),*

$$\liminf_{n\to\infty}\frac{1}{b_n}\log\mathbb{E}\exp\bigg\{\frac{b_n\log n}{n}\int_{\mathbb{R}^2}f\bigg(\sqrt{\frac{b_n}{n}}x\bigg)\mathbb{1}_{\{[x]\in S[0,n]\}}\,dx\bigg\}\ge\sup_{g\in\mathcal{F}_2}\bigg\{2\pi\sqrt{\det(\Gamma)}\int_{\mathbb{R}^2}f(x)g^2(x)\,dx-\frac12\int_{\mathbb{R}^2}\langle\nabla g(x),\Gamma\nabla g(x)\rangle\,dx\bigg\}. \tag{4.24}$$

(ii) *As* $d=3$, $f\in\mathcal{L}^2(\mathbb{R}^3)$ *and* $\{b_n\}$ *satisfies* (1.12),

$$\liminf_{n\to\infty}\frac{1}{b_n}\log\mathbb{E}\exp\left\{\frac{b_n}{n}\int_{\mathbb{R}^3}f\left(\sqrt{\frac{b_n}{n}}x\right)\mathbb{1}_{\{[x]\in S[0,n]\}}\,dx\right\}$$

(4.25)

$$\ge\sup_{g\in\mathcal{F}_3}\left\{\gamma(S)\int_{\mathbb{R}^3}f(x)g^2(x)\,dx-\frac{1}{2}\int_{\mathbb{R}^3}\langle\nabla g(x),\Gamma\nabla g(x)\rangle\,dx\right\}.$$

PROOF. Due to similarity we only prove (4.24). We first assume that $\{S(n)\}$ is aperiodic. By uniform continuity

$$\left|\int_{\mathbb{R}^2}f\left(\sqrt{\frac{b_n}{n}}x\right)\mathbb{1}_{\{[x]\in S[0,n]\}}\,dx-\sum_{x\in S[0,n]}f\left(\sqrt{\frac{b_n}{n}}x\right)\right|\le\theta_n\#\{S[0,n]\},\tag{4.26}$$

where $\{\theta_n\}$ is a deterministic positive sequence with $\theta_n\to 0$ as $n\to\infty$. Recall that $t_n$ and $\Delta_i$ are defined by (4.16). Notice that

$$\mathbb{E}\exp\left\{\theta\frac{b_n\log n}{n}\#\{S[0,n]\}\right\}\le\left(\mathbb{E}\exp\left\{\theta\frac{b_n\log n}{n}\#\{S[0,t_n]\}\right\}\right)^{b_n+1}.$$

By Lemma 3,

$$\limsup_{n\to\infty}\frac{1}{b_n}\log\mathbb{E}\exp\left\{\theta\frac{b_n\log n}{n}\#\{S[0,n]\}\right\}\le\Lambda(\theta)\qquad(\theta>0)\tag{4.27}$$

where $\Lambda(\theta)\to 0$ as $\theta\to 0^+$. By (4.26), (4.27) and a standard argument of exponential approximation, (4.24) is equivalent to

$$\liminf_{n\to\infty}\frac{1}{b_n}\log\mathbb{E}\exp\left\{\frac{b_n\log n}{n}\sum_{x\in S[0,n]}f\left(\sqrt{\frac{b_n}{n}}x\right)\right\}$$

(4.28)

$$\ge\sup_{g\in\mathcal{F}_2}\left\{2\pi\sqrt{\det(\Gamma)}\int_{\mathbb{R}^2}f(x)g^2(x)\,dx-\frac{1}{2}\int_{\mathbb{R}^2}\langle\nabla g(x),\Gamma\nabla g(x)\rangle\,dx\right\}.$$

To prove (4.28), notice that

$$\sum_{x\in S[0,n]}f\left(\sqrt{\frac{b_n}{n}}x\right)$$

$$\ge\sum_{i=1}^{[b_n]}\sum_{x\in S(\Delta_j)}f\left(\sqrt{\frac{b_n}{n}}x\right)-\sum_{1\le j<k\le[b_n]}\sum_{x\in S(\Delta_j)\cap S(\Delta_k)}f\left(\sqrt{\frac{b_n}{n}}x\right).$$

Therefore, for any given $\varepsilon > 0$,

$$
\begin{aligned}
&\mathbb{E}\exp\left\{\frac{b_n \log n}{n}\sum_{x\in S[0,n]} f\left(\sqrt{\frac{b_n}{n}}x\right)\right\}\\
&\quad \ge \mathbb{E}\exp\left\{\frac{b_n\log n}{n}\left(\sum_{i=1}^{[b_n]}\sum_{x\in S(\Delta_i)} f\left(\sqrt{\frac{b_n}{n}}x\right)\right.\right.\\
&\qquad\qquad \left.\left. - \sum_{1\le j<k\le [b_n]}\sum_{x\in S(\Delta_j)\cap S(\Delta_k)} f\left(\sqrt{\frac{b_n}{n}}x\right)\right)\right\}\\
&\quad \ge e^{-\varepsilon b_n}\mathbb{E}\exp\left\{\frac{b_n\log n}{n}\sum_{i=1}^{[b_n]}\sum_{x\in S(\Delta_i)} f\left(\sqrt{\frac{b_n}{n}}x\right)\right\}\\
&\qquad - \mathbb{E}\left[\exp\left\{\frac{b_n\log n}{n}\sum_{i=1}^{[b_n]}\sum_{x\in S(\Delta_i)} f\left(\sqrt{\frac{b_n}{n}}x\right)\right\};\right.\\
&\qquad\qquad \left.\sum_{1\le j<k\le [b_n]}\sum_{x\in S(\Delta_j)\cap S(\Delta_k)} f\left(\sqrt{\frac{b_n}{n}}x\right)\ge \varepsilon\frac{n}{\log n}\right]\\
&\quad = (I) - (II) \qquad \text{(say).}
\end{aligned}
\tag{4.29}
$$

By Lemma 5,

$$
\begin{aligned}
\liminf_{n\to\infty}\frac{1}{b_n}\log(I) \ge -\varepsilon + \sup_{g\in\mathcal{F}_2}&\left\{2\pi\sqrt{\det(\Gamma)}\int_{\mathbb{R}^2} f(x)g^2(x)\,dx\right.\\
&\left. - \frac{1}{2}\int_{\mathbb{R}^2}\langle\nabla g(x),\Gamma\nabla g(x)\rangle\,dx\right\}.
\end{aligned}
\tag{4.30}
$$

By the Cauchy–Schwarz inequality,

$$
\begin{aligned}
(II) \le &\left[\mathbb{E}\exp\left\{\frac{2b_n\log n}{n}\sum_{i=1}^{[b_n]}\sum_{x\in S(\Delta_i)} f\left(\sqrt{\frac{b_n}{n}}x\right)\right\}\right]^{1/2}\\
&\times\left[\mathbb{P}\left\{\sum_{1\le j<k\le [b_n]}\sum_{x\in S(\Delta_j)\cap S(\Delta_k)} f\left(\sqrt{\frac{b_n}{n}}x\right)\ge\varepsilon\frac{n}{\log n}\right\}\right]^{1/2}.
\end{aligned}
\tag{4.31}
$$

Notice that

$$
\sum_{1\le j<k\le [b_n]}\sum_{x\in S(\Delta_j)\cap S(\Delta_k)} f\left(\sqrt{\frac{b_n}{n}}x\right) \le \|f\|_\infty \sum_{1\le j<k\le [b_n]} \#\{S(\Delta_j)\cap S(\Delta_k)\}.
$$

By Lemma 6,

$$\limsup_{n\to\infty}\frac{1}{b_n}\log\mathbb{P}\Big\{\sum_{1\le j<k\le[b_n]}\sum_{x\in S(\Delta_j)\cap S(\Delta_k)} f\Big(\sqrt{\frac{b_n}{n}}x\Big)\ge\varepsilon\frac{n}{\log n}\Big\}=-\infty. \tag{4.32}$$

In view of (4.29)–(4.32), it remains to prove

$$\limsup_{n\to\infty}\frac{1}{b_n}\log\mathbb{E}\exp\left\{\frac{2b_n\log n}{n}\sum_{i=1}^{[b_n]}\sum_{x\in S(\Delta_i)} f\Big(\sqrt{\frac{b_n}{n}}x\Big)\right\}<\infty.$$

By the exponential approximation used earlier, this is equivalent to

$$\limsup_{n\to\infty}\frac{1}{b_n}\log\mathbb{E}\exp\left\{\frac{2b_n\log n}{n}\int_{\mathbb{R}^2} f\Big(\sqrt{\frac{b_n}{n}}x\Big)\sum_{i=1}^{[b_n]}\mathbb{1}_{\{[x]\in S(\Delta_i)\}}\,dx\right\}<\infty. \tag{4.33}$$

For any integer $m\ge 1$,

$$\begin{aligned}
&\mathbb{E}\left(\int_{\mathbb{R}^2} f\Big(\sqrt{\frac{b_n}{n}}x\Big)\sum_{i=1}^{[b_n]}\mathbb{1}_{\{[x]\in S(\Delta_i)\}}\,dx\right)^m\\
&=\Big(\frac{n}{b_n}\Big)^m\mathbb{E}\left(\int_{\mathbb{R}^2} f(x)\sum_{i=1}^{[b_n]}\mathbb{1}_{\{[\sqrt{n/b_n}x]\in S(\Delta_i)\}}\,dx\right)^m\\
&=\Big(\frac{n}{b_n}\Big)^m\int_{(\mathbb{R}^2)^m}dx_1\cdots dx_m\Big(\prod_{k=1}^m f(x_k)\Big)\mathbb{E}\prod_{k=1}^m\sum_{i=1}^{[b_n]}\mathbb{1}_{\{[\sqrt{n/b_n}x_k]\in S(\Delta_i)\}}\\
&\le\|f\|_q^m\Big(\frac{n}{b_n}\Big)^m\left\{\int_{(\mathbb{R}^2)^m}dx_1\cdots dx_m\left(\mathbb{E}\prod_{k=1}^m\sum_{i=1}^{[b_n]}\mathbb{1}_{\{[\sqrt{n/b_n}x_k]\in S(\Delta_i)\}}\right)^p\right\}^{1/p} \qquad (4.34)\\
&=\|f\|_q^m\Big(\frac{n}{b_n}\Big)^{(p-1)/pm}\left\{\int_{(\mathbb{R}^2)^m}dx_1\cdots dx_m\left(\mathbb{E}\prod_{k=1}^m\sum_{i=1}^{[b_n]}\mathbb{1}_{\{[x_k]\in S(\Delta_i)\}}\right)^p\right\}^{1/p}\\
&=\|f\|_q^m\Big(\frac{n}{b_n}\Big)^{(p-1)/pm}\left\{\sum_{x_1,\dots,x_m}\left(\mathbb{E}\prod_{k=1}^m\sum_{i=1}^{[b_n]}\mathbb{1}_{\{x_k\in S(\Delta_i)\}}\right)^p\right\}^{1/p}\\
&=\|f\|_q^m\Big(\frac{n}{b_n}\Big)^{(p-1)/pm}\left\{\mathbb{E}\left(\sum_{x\in\mathbb{Z}^2}\prod_{j=1}^p\sum_{i=1}^{[b_n]}\mathbb{1}_{\{x\in S_j(\Delta_i)\}}\right)^m\right\}^{1/p}.
\end{aligned}$$

Similar to (4.23),

$$\limsup_{n\to\infty}\frac{1}{b_n}\log\sum_{m=0}^{\infty}\frac{(2\|f\|_q)^m}{m!}\left(\frac{b_n(\log n)^p}{n}\right)^{m/p}$$
$$\times\left\{\mathbb{E}\left(\sum_{x\in\mathbb{Z}^2}\prod_{j=1}^{p}\sum_{i=1}^{[b_n]}\mathbb{1}_{\{x\in S_j(\Delta_i)\}}\right)^m\right\}^{1/p}<\infty.$$

So (4.33) follows from (4.34).

We now prove (4.24) without assuming aperiodicity. Let $0<\eta<1$ be fixed and let $\{\delta_n\}_{n\ge1}$ be i.i.d. Bernoulli random variables with the common law:

$$\mathbb{P}\{\delta_1=0\}=1-\mathbb{P}\{\delta_1=1\}=\eta.$$

We assume independence between $\{S(n)\}$ and $\{\delta_n\}$.

Define the renewal sequence $\{\sigma_k\}_{k\ge0}$ by

$$\sigma_0=0\quad\text{and}\quad\sigma_{k+1}=\inf\{n>\sigma_k;\delta_n=1\}.$$

Then $\{\sigma_k-\sigma_{k-1}\}_{k\ge1}$ is an i.i.d. sequence with common distribution

$$\mathbb{P}\{\sigma_1=n\}=(1-\eta)\eta^{n-1},\qquad n=1,2,\dots.$$

Consider the random walk $\tilde S(n)=S(\sigma_n)$. $\{\tilde S(n)\}$ is symmetric with covariance

$$\mathrm{Cov}\,(S(\sigma_1),S(\sigma_1))=(\mathbb{E}\sigma_1)\Gamma=(1-\eta)^{-1}\Gamma.$$

By the fact that

$$\mathbb{P}\{S(\sigma_1)=0\}=(1-\eta)\sum_{k=1}^{\infty}\eta^{k-1}\mathbb{P}\{S(k)=0\}>0,$$

$\{\tilde S(n)\}$ is aperiodic. Applying what we have proved to $\{\tilde S(n)\}$,

$$\liminf_{n\to\infty}\frac{1}{b_n}\log\mathbb{E}\exp\left\{\frac{b_n\log n}{n}\int_{\mathbb{R}^2}f\left(\sqrt{\frac{b_n}{n}}x\right)\mathbb{1}_{\{[x]\in\tilde S[0,n]\}}\,dx\right\}$$
$$\ge\sup_{g\in\mathcal{F}_2}\left\{2\pi\frac{\sqrt{\det(\Gamma)}}{1-\eta}\int_{\mathbb{R}^2}f(x)g^2(x)\,dx\right.$$
$$\left.-\frac{1}{2(1-\eta)}\int_{\mathbb{R}^2}\langle\nabla g(x),\Gamma\nabla g(x)\rangle\,dx\right\}$$
$$\ge\sup_{g\in\mathcal{F}_2}\left\{2\pi\sqrt{\det(\Gamma)}\int_{\mathbb{R}^2}f(x)g^2(x)\,dx-\frac12\int_{\mathbb{R}^2}\langle\nabla g(x),\Gamma\nabla g(x)\rangle\,dx\right\}.$$

Notice that

$$\tilde S[0,n]=\{S(\sigma_0),\dots,S(\sigma_n)\}\subset S[0,\sigma_n].$$

Given $\varepsilon > 0$,

$$\begin{aligned}
&\mathbb{E}\exp\left\{\frac{b_n \log n}{n}\int_{\mathbb{R}^2} f\left(\sqrt{\frac{b_n}{n}}x\right)\mathbb{1}_{\{[x]\in S[0,[(1+\varepsilon)n]]\}}\,dx\right\}\\
&\quad\ge \mathbb{E}\exp\left\{\frac{b_n \log n}{n}\int_{\mathbb{R}^2} f\left(\sqrt{\frac{b_n}{n}}x\right)\mathbb{1}_{\{[x]\in \tilde{S}[0,n]\}}\,dx\right\}\\
&\qquad - \mathbb{E}\exp\left\{\frac{b_n \log n}{n}\int_{\mathbb{R}^2} f\left(\sqrt{\frac{b_n}{n}}x\right)\mathbb{1}_{\{[x]\in \tilde{S}[0,n]\}}\,dx\right\}\mathbb{1}_{\{\sigma_n\ge(1+\varepsilon)n\}}.
\end{aligned}$$

By Cramér large deviation [Theorem 2.2.3 of Dembo and Zeitouni (1998)] as $(1-\eta)^{-1} < 1+\varepsilon$ there is $u > 0$ such that

$$\mathbb{P}\{\sigma_n \ge (1+\varepsilon)n\} \le e^{-un}$$

for sufficiently large $n$. By (4.33) [with $S(n)$ being replaced by $\tilde{S}(n)$] and the Cauchy–Schwarz inequality, therefore,

$$\limsup_{n\to\infty}\frac{1}{b_n}\log\mathbb{E}\exp\left\{\frac{b_n \log n}{n}\int_{\mathbb{R}^2} f\left(\sqrt{\frac{b_n}{n}}x\right)\mathbb{1}_{\{[x]\in \tilde{S}[0,n]\}}\,dx\right\}\mathbb{1}_{\{\sigma_n\ge(1+\varepsilon)n\}} = -\infty.$$

Hence,

$$\begin{aligned}
&\liminf_{n\to\infty}\frac{1}{b_n}\log\mathbb{E}\exp\left\{\frac{b_n \log n}{n}\int_{\mathbb{R}^2} f\left(\sqrt{\frac{b_n}{n}}x\right)\mathbb{1}_{\{[x]\in S[0,[(1+\varepsilon)n]]\}}\,dx\right\}\\
&\quad\ge \sup_{g\in\mathcal{F}_2}\left\{2\pi\sqrt{\det(\Gamma)}\int_{\mathbb{R}^2} f(x)g^2(x)\,dx - \frac{1}{2}\int_{\mathbb{R}^2}\langle\nabla g(x),\Gamma\nabla g(x)\rangle\,dx\right\}.
\end{aligned}$$

Replacing $[(1+\varepsilon)n]$ by $n$ and $f(x)$ by $(1+\varepsilon)^{-1}f((1+\varepsilon)^{-1/2}x)$, we have

$$\begin{aligned}
&\liminf_{n\to\infty}\frac{1}{b_n}\log\mathbb{E}\exp\left\{\frac{b_n \log n}{n}\int_{\mathbb{R}^2} f\left(\sqrt{\frac{b_n}{n}}x\right)\mathbb{1}_{\{[x]\in S[0,n]\}}\,dx\right\}\\
&\quad\ge \sup_{g\in\mathcal{F}_2}\left\{2\pi\sqrt{\det(\Gamma)}(1+\varepsilon)^{-1}\int_{\mathbb{R}^2} f(x)g^2(x)\,dx\right.\\
&\qquad\qquad\left. - \frac{1}{2}\int_{\mathbb{R}^2}\langle\nabla g(x),\Gamma\nabla g(x)\rangle\,dx\right\}.
\end{aligned}$$

Letting $\varepsilon\to 0^+$ gives (4.21). □

We are finally ready to prove (4.1) and (4.2). Due to similarity we only prove (4.1). Notice that

$$J_n = \sum_{x\in\mathbb{Z}^2}\prod_{j=1}^{p}\mathbb{1}_{\{x\in S_j[0,n]\}}.$$

For any nonnegative, bounded and uniformly continuous function $f$ on $\mathbb{R}^2$ with $\|f\|_q = 1$, a procedure similar to (4.34) gives

$$\left(\frac{n}{b_n}\right)^{(p-1)/pm} (\mathbb{E}J_n^m)^{1/p} \geq \mathbb{E}\left(\int_{\mathbb{R}^2} f\left(\sqrt{\frac{b_n}{n}}x\right) \mathbb{1}_{\{[x]\in S[0,n]\}}\, dx\right)^m, \qquad m = 0, 1, \ldots.$$

Therefore,

$$\sum_{m=0}^{\infty} \frac{\theta^m}{m!}\left(\frac{b_n \log^p n}{n}\right)^{m/p} (\mathbb{E}J_n^m)^{1/p} \geq \mathbb{E}\exp\left\{\theta \frac{b_n \log n}{n}\int_{\mathbb{R}^2} f\left(\sqrt{\frac{b_n}{n}}x\right)\mathbb{1}_{\{[x]\in S[0,n]\}}\, dx\right\}.$$

By Theorem 8,

$$\liminf_{n\to\infty} \frac{1}{b_n} \log \sum_{m=0}^{\infty} \frac{\theta^m}{m!}\left(\frac{b_n \log^p n}{n}\right)^{m/p} (\mathbb{E}J_n^m)^{1/p} \geq \sup_{g\in\mathcal{F}_2}\left\{2\pi\theta\sqrt{\det(\Gamma)}\int_{\mathbb{R}^2} f(x) g^2(x)\, dx - \frac{1}{2}\int_{\mathbb{R}^2}\langle \nabla g(x), \Gamma \nabla g(x)\rangle\, dx\right\}.$$

Taking the supremum over all nonnegative, bounded and uniformly continuous functions $f$ on $\mathbb{R}^2$ with $\|f\|_q = 1$, the right-hand side becomes

$$\begin{aligned} &\sup_{g\in\mathcal{F}_2}\left\{2\pi\theta\sqrt{\det(\Gamma)}\left(\int_{\mathbb{R}^2}|g(x)|^{2p}\, dx\right)^{1/p} - \frac{1}{2}\int_{\mathbb{R}^2}\langle\nabla g(x), \Gamma\nabla g(x)\rangle\, dx\right\} \\ &\quad = (2\pi\theta)^p\sqrt{\det(\Gamma)}\sup_{g\in\mathcal{F}_2}\left\{\left(\int_{\mathbb{R}^2}|h(x)|^{2p}\, dx\right)^{1/p} - \frac{1}{2}\int_{\mathbb{R}^2}|\nabla h(x)|^2\, dx\right\} \qquad (4.35) \\ &\quad = \frac{1}{p}\left(\frac{2(p-1)}{p}\right)^{p-1}(2\pi\theta)^p\sqrt{\det(\Gamma)}\kappa(2,p)^{2p}, \end{aligned}$$

where the first equality follows from the substitution $g(x) = \sqrt{|\det A|}h(Ax)$ with the $2\times 2$ matrix $A$ satisfying

$$A^\tau \Gamma A = (2\pi\theta)^p\sqrt{\det(\Gamma)}\mathbf{I}_2$$

with $\mathbf{I}_2$ being the $2\times 2$ identity matrix, and the second equality follows from Lemma A.2 in Chen (2004).

**5. Law of the iterated logarithm.** We prove Theorem 3 in this section. With the moderate deviations given in Theorems 1 and 2, the proof of the upper bound is just a standard practice of the Borel–Cantelli lemma. So we only give proof to the lower bounds. That is, we prove:

As $d=2$ and $p\ge 2$,

$$\limsup_{n\to\infty}\frac{(\log n)^p}{n(\log\log n)^{p-1}}J_n\ge (2\pi)^p\left(\frac{2}{p}\right)^{p-1}\sqrt{\det(\Gamma)}\kappa(2,p)^{2p}\qquad\text{a.s.}\tag{5.1}$$

As $d=3$ and $p=2$,

$$\limsup_{n\to\infty}\frac{1}{\sqrt{n(\log\log n)^3}}J_n\ge\gamma(S)^2\det(\Gamma)^{-1/2}\kappa(3,2)^4\qquad\text{a.s.}\tag{5.2}$$

By the technology used in the proof of Theorem 8, which extends the lower bound established under aperiodicity to the general case, we may assume aperiodicity in the proof given below.

For given $\bar{x}=(x_1,\dots,x_p)\in(\mathbb{Z}^d)^p$, we introduce the notation $\mathbb{P}^{\bar{x}}$ for the probability induced by the random walks $S_1(n),\dots,S_p(n)$ in the case when $S_1(n),\dots,S_p(n)$ start at $x_1,\dots,x_p$, respectively. The notation $\mathbb{E}^{\bar{x}}$ denotes the expectation correspondent to $\mathbb{P}^{\bar{x}}$. To be consistent with the notation we used before, we have $\mathbb{P}^{(0,\dots,0)}=\mathbb{P}$ and $\mathbb{E}^{(0,\dots,0)}=\mathbb{E}$. Write

$$\|\bar{x}\|=\max_{1\le j\le p}|x_j|,\qquad \bar{x}=(x_1,\dots,x_p)\in(\mathbb{R}^d)^p.$$

LEMMA 7. *Under the conditions in Theorem 1,*

$$\begin{aligned}\liminf_{n\to\infty}\frac{1}{b_n}\log\inf_{\|\bar{x}\|\le\sqrt{n/b_n}}\mathbb{P}^{\bar{x}}\left\{J_n\ge\lambda\frac{n}{(\log n)^p}b_n^{p-1}\right\}&\\ \ge-\frac{p}{2}(2\pi)^{-p/(p-1)}&\\ \times\det(\Gamma)^{-1/(2(p-1))}\kappa(2,p)^{-2p/(p-1)}\lambda^{1/(p-1)}&\qquad(\lambda>0).\end{aligned}\tag{5.3}$$

*Under the conditions in Theorem 2,*

$$\begin{aligned}\liminf_{n\to\infty}\frac{1}{b_n}\log\inf_{\|\bar{x}\|\le\sqrt{n/b_n}}\mathbb{P}^{\bar{x}}\{J_n\ge\lambda\sqrt{nb_n^3}\}&\\ \ge-\det(\Gamma)^{1/3}\gamma(S)^{-4/3}\kappa(3,2)^{-8/3}\lambda^{2/3}&\qquad(\lambda>0).\end{aligned}\tag{5.4}$$

PROOF. Due to similarity we only prove (5.3). For given $\bar{y}=(y_1,\dots,y_p)\in(\mathbb{Z}^2)^p$ and $m,n\ge 1$,

$$\begin{aligned}\mathbb{E}^{\bar{y}}J_n^m&=\sum_{x_1,\dots,x_m}\prod_{j=1}^p\mathbb{E}\prod_{k=1}^m\mathbb{1}_{\{y_j+x_k\in S[0,n]\}}\\ &\le\prod_{j=1}^p\left(\sum_{x_1,\dots,x_m}\left[\mathbb{E}\prod_{k=1}^m\mathbb{1}_{\{y_j+x_k\in S[0,n]\}}\right]^p\right)^{1/p}\\ &=\sum_{x_1,\dots,x_m}\left[\mathbb{E}\prod_{k=1}^m\mathbb{1}_{\{x_k\in S[0,n]\}}\right]^p=\mathbb{E}J_n^m.\end{aligned}$$

By (3.1) we have

$$\limsup_{n\to\infty}\frac{1}{b_n}\log\sum_{m=0}^{\infty}\frac{\theta^m}{m!}\left(\frac{b_n\log^p n}{n}\right)^{m/p}\left(\sup_{\bar{y}}\mathbb{E}^{\bar{y}}J_n^m\right)^{1/p}$$
$$\le\frac{1}{p}\left(\frac{2(p-1)}{p}\right)^{p-1}(2\pi\theta)^p\sqrt{\det(\Gamma)}\kappa(2,p)^{2p}\qquad(\theta>0).$$

It is easy to see from Theorem 4 that we will have (5.3) if we can prove

$$\liminf_{n\to\infty}\frac{1}{b_n}\log\sum_{m=0}^{\infty}\frac{\theta^m}{m!}\left(\frac{b_n\log^p n}{n}\right)^{m/p}\left(\inf_{\|\bar{y}\|\le\sqrt{n/b_n}}\mathbb{E}^{\bar{y}}J_n^m\right)^{1/p}$$
$$\ge\frac{1}{p}\left(\frac{2(p-1)}{p}\right)^{p-1}(2\pi\theta)^p\sqrt{\det(\Gamma)}\kappa(2,p)^{2p}\tag{5.5}$$

for every $\theta>0$.

Let $\varepsilon>0$ be fixed for a moment. For any sets $A,B\subset\mathbb{Z}^2$, $A+B$ is defined as the set $\{x+y;x\in A \text{ and } y\in B\}$. In particular, $x+B\equiv\{x\}+B$ for any $x\in\mathbb{Z}^2$. Write

$$B_n(x)=\{y;|y-x|\le\varepsilon\sqrt{n/b_n}\},\qquad x\in\mathbb{Z}^2,$$

and set $B_n=B_n(0)$.

For any function $f$ on $\mathbb{R}^2$, write

$$f_\varepsilon(x)=\frac{1}{\pi\varepsilon^2}\int_{\{|y|\le\varepsilon\}}f(x+y)\,dy$$

whenever the integral on the right-hand side makes sense.

Define

$$J_n(\varepsilon)=\sum_{x\in\mathbb{Z}^2}\prod_{j=1}^{p}\left(\frac{1}{\#(B_n)}\mathbb{1}_{\{x\in S_j[0,n]+B_n\}}\right).$$

Let $f$ be a nonnegative, bounded and uniformly continuous function on $\mathbb{R}^2$ with $\|f\|_q=1$:

$$\int_{\mathbb{R}^2}f\left(\sqrt{\frac{b_n}{n}}x\right)\frac{1}{\#(B_n)}\mathbb{1}_{\{[x]\in S[0,n]+B_n\}}\,dx$$
$$=\int_{\mathbb{R}^2}f\left(\sqrt{\frac{b_n}{n}}x\right)\frac{1}{\#(B_n)}\sum_{y\in B_n}\mathbb{1}_{\{[x-y]\in S[0,n]\}}\,dx$$
$$=\int_{\mathbb{R}^2}\mathbb{1}_{\{[x]\in S[0,n]\}}\left(\frac{1}{\#(B_n)}\sum_{y\in B_n}f\left(\sqrt{\frac{b_n}{n}}(x+y)\right)\right)dx$$

$$= (1+o(1))\int_{\mathbb{R}^2} \mathbb{1}_{\{[x]\in S[0,n]\}} f_\varepsilon\left(\sqrt{\frac{b_n}{n}}x\right) dx$$

where $o(1)$ is bounded by a deterministic sequence that approaches to zero as $n \to \infty$.

Similar to (4.34), for any integer $m \ge 1$,

$$\begin{aligned}
&\left(\frac{n}{b_n}\right)^{(p-1)/pm} (\mathbb{E}J_n(\varepsilon)^m)^{1/p} \\
&\qquad \ge \mathbb{E}\left(\int_{\mathbb{R}^2} f\left(\sqrt{\frac{b_n}{n}}x\right) \frac{1}{\#(B_n)} \mathbb{1}_{\{[x]\in S[0,n]+B_n\}} dx\right)^m \\
&\qquad \ge (1+o(1))^m \mathbb{E}\left(\int_{\mathbb{R}^2} \mathbb{1}_{\{[x]\in S[0,n]\}} f_\varepsilon\left(\sqrt{\frac{b_n}{n}}x\right) dx\right)^m.
\end{aligned}$$

Therefore,

$$\begin{aligned}
&\sum_{m=0}^\infty \frac{\theta^m}{m!}\left(\frac{b_n \log^p n}{n}\right)^{m/p} (\mathbb{E}J_n(\varepsilon)^m)^{1/p} \\
&\qquad \ge \mathbb{E}\exp\left\{(1+o(1))\theta \frac{b_n \log n}{n} \int_{\mathbb{R}^2} f_\varepsilon\left(\sqrt{\frac{b_n}{n}}x\right) \mathbb{1}_{\{[x]\in S[0,n]\}} dx\right\}.
\end{aligned}$$

By Theorem 8,

$$\begin{aligned}
&\liminf_{n\to\infty} \frac{1}{b_n} \log \sum_{m=0}^\infty \frac{\theta^m}{m!}\left(\frac{b_n \log^p n}{n}\right)^{m/p} (\mathbb{E}J_n(\varepsilon)^m)^{1/p} \\
&\quad \ge \sup_{g\in\mathcal{F}_2}\left\{2\pi\theta\sqrt{\det(\Gamma)}\int_{\mathbb{R}^2} f_\varepsilon(x) g^2(x)\, dx - \frac{1}{2}\int_{\mathbb{R}^2}\langle\nabla g(x), \Gamma\nabla g(x)\rangle\, dx\right\} \\
&\quad = \sup_{g\in\mathcal{F}_2}\left\{2\pi\theta\sqrt{\det(\Gamma)}\int_{\mathbb{R}^2} f(x) (g^2)_\varepsilon(x)\, dx - \frac{1}{2}\int_{\mathbb{R}^2}\langle\nabla g(x), \Gamma\nabla g(x)\rangle\, dx\right\}.
\end{aligned}$$

Taking the supremum over all nonnegative, bounded and uniformly continuous functions $f$ on $\mathbb{R}^2$ with $\|f\|_q = 1$ gives

$$\begin{aligned}
&\liminf_{n\to\infty} \frac{1}{b_n} \log \sum_{m=0}^\infty \frac{\theta^m}{m!}\left(\frac{b_n \log^p n}{n}\right)^{m/p} (\mathbb{E}J_n(\varepsilon)^m)^{1/p} \\
&\quad \ge \sup_{g\in\mathcal{F}_2}\left\{2\pi\theta\sqrt{\det(\Gamma)}\left(\int_{\mathbb{R}^2} |(g^2)_\varepsilon(x)|^p\, dx\right)^{1/p} \right. \\
&\qquad\qquad \left. - \frac{1}{2}\int_{\mathbb{R}^2}\langle\nabla g(x), \Gamma\nabla g(x)\rangle\, dx\right\}.
\end{aligned} \tag{5.6}$$

Take $t_n = [n/b_n]$. To prove (5.5), notice that

$$\begin{aligned}
\mathbb{E}^{\bar{y}} J_n^m &\geq \mathbb{E}\bigg(\sum_x \prod_{j=1}^p \mathbb{1}_{\{y_j + x \in S_j[t_n, n]\}}\bigg)^m \\
&= \sum_{x_1,\ldots,x_m} \prod_{j=1}^p \mathbb{E}\bigg(\prod_{k=1}^m \mathbb{1}_{\{y_j + x_k \in S[t_n,n]\}}\bigg) \\
&\geq \sum_{x_1,\ldots,x_m} \prod_{j=1}^p \mathbb{E}\bigg(\prod_{k=1}^m \sum_{z\in B_n(y_j)} \mathbb{1}_{\{S(t_n)=z\}} \cdot \mathbb{1}_{\{y_j - z + x_k \in S'[0,n-t_n]\}}\bigg)
\end{aligned}$$

where $S'(k) = S(k+t_n) - S(t_n)$ $(k = 1, 2, \ldots)$. By the identity,

$$\begin{aligned}
&\prod_{k=1}^m \sum_{z\in B_n(y_j)} \mathbb{1}_{\{S(t_n)=z\}} \cdot \mathbb{1}_{\{y_j - y + x_k \in S'[0,n-t_n]\}} \\
&\qquad = \sum_{y\in B_n(y_j)} \mathbb{1}_{\{S(t_n)=z\}} \cdot \prod_{k=1}^m \mathbb{1}_{\{y_j - z + x_k \in S'[0,n-t_n]\}}
\end{aligned}$$

and therefore by independence,

$$\begin{aligned}
&\mathbb{E}\bigg(\prod_{k=1}^m \sum_{z\in B_n(y_j)} \mathbb{1}_{\{S(t_n)=z\}} \cdot \mathbb{1}_{\{y_j - z + x_k \in S'[0,n-t_n]\}}\bigg) \\
&\quad = \sum_{z\in B_n(y_j)} \mathbb{P}\{S(t_n) = z\} \cdot \mathbb{E}\bigg[\prod_{k=1}^m \mathbb{1}_{\{y_j - z + x_k \in S[0,n-t_n]\}}\bigg] \\
&\quad \geq \min_{1\leq j\leq p} \inf_{z\in B_n(y_j)} \{\mathbb{P}\{S(t_n) = z\}\} \sum_{z\in B_n} \mathbb{E}\bigg[\prod_{k=1}^m \mathbb{1}_{\{x_k - z \in S[0,n-t_n]\}}\bigg] \\
&\quad = \gamma_n \sum_{z\in B_n} \mathbb{E}\bigg[\prod_{k=1}^m \mathbb{1}_{\{x_k - z \in S[0,n-t_n]\}}\bigg] \qquad \text{(say)}.
\end{aligned}$$

Hence,

$$\begin{aligned}
\mathbb{E}^{\bar{y}} J_n^m &\geq \gamma_n^p \sum_{x_1,\ldots,x_m} \bigg(\sum_{z\in B_n} \mathbb{E}\bigg[\prod_{k=1}^m \mathbb{1}_{\{x_k - z \in S[0,n-t_n]\}}\bigg]\bigg)^p \\
&= \sum_{x_1,\ldots,x_m} \sum_{z_1,\ldots,z_p\in B_n} \prod_{j=1}^p \mathbb{E}\bigg[\prod_{k=1}^m \mathbb{1}_{\{x_k - z_j \in S[0,n-t_n]\}}\bigg] \\
&= \gamma_n^p \sum_{z_1,\ldots,z_p\in B_n} \sum_{x_1,\ldots,x_m} \mathbb{E}\bigg[\prod_{k=1}^m \prod_{j=1}^p \mathbb{1}_{\{x_k - z_j \in S_j[0,n-t_n]\}}\bigg]
\end{aligned}$$

$$= \gamma_n^p \sum_{z_1,\dots,z_p\in B_n} \mathbb{E}\biggl(\sum_x \prod_{j=1}^p \mathbb{1}_{\{x-z_j\in S_j[0,n-t_n]\}}\biggr)^m$$

$$\ge \gamma_n^p \#\{B_n\}^p \mathbb{E}\biggl(\frac{1}{\#\{B_n\}^p}\sum_{z_1,\dots,z_p\in B_n}\sum_x \prod_{j=1}^p \mathbb{1}_{\{x-z_j\in S_j[0,n-t_n]\}}\biggr)^m$$

$$= \gamma_n^p \#\{B_n\}^p \mathbb{E}\biggl(\sum_x \prod_{j=1}^p \frac{1}{\#\{B_n\}}\mathbb{1}_{\{x\in S_j[0,n-t_n]+B_n\}}\biggr)^m,$$

where the fifth step follows from Jensen's inequality. By (4.3) (with $d=2$ and $n$ replaced by $t_n$),

$$\gamma_n = \frac{1}{t_n}\min_{1\le j\le p}\inf_{|y_j|\le\sqrt{n/b_n}}\inf_{z\in B_n(y_j)}\biggl[\exp\biggl\{-\frac{1}{2t_n}\langle y,\Gamma^{-1}y\rangle\biggr\}+o(1)\biggr]\ge ct_n^{-1}.$$

We have proved that there is a $\delta=\delta(\varepsilon)>0$, such that for any integer $m\ge 0$ and $n\ge 1$,

$$\inf_{\|\bar y\|\le\sqrt{n/b_n}}\mathbb{E}^{\bar y}J_n^m \ge \delta\mathbb{E}\biggl(\sum_x\prod_{j=1}^p\frac{1}{\#\{B_n\}}\mathbb{1}_{\{x\in S_j[0,n-t_n]+B_n\}}\biggr)^m.$$

By (5.6) (with $n$ replaced by $n-t_n$),

$$\liminf_{n\to\infty}\frac{1}{b_n}\log\sum_{m=0}^\infty\frac{\theta^m}{m!}\biggl(\frac{b_n\log^p n}{n}\biggr)^{m/p}\Bigl(\inf_{\|\bar y\|\le\sqrt{n/b_n}}\mathbb{E}^{\bar y}J_n^m\Bigr)^{1/p}$$

$$\ge \sup_{g\in\mathcal{F}_2}\biggl\{2\pi\theta\sqrt{\det(\Gamma)}\biggl(\int_{\mathbb{R}^2}|(g^2)_\varepsilon(x)|^p\,dx\biggr)^{1/p} - \frac{1}{2}\int_{\mathbb{R}^2}\langle\nabla g(x),\Gamma\nabla g(x)\rangle\,dx\biggr\}.$$

Finally, we let $\varepsilon\to 0^+$ on the right-hand side. Then (5.5) follows from (4.35). $\square$

We only prove (5.1) as the proof of (5.2) is analogous. Let $n_k=k^k$. We first show that for any $\lambda<(2\pi)^p(\frac{2}{p})^{p-1}\sqrt{\det(\Gamma)}\kappa(2,p)^{2p}$,

$$\limsup_{k\to\infty}\frac{(\log n_{k+1})^p}{n_{k+1}\log\log n_{k+1}}\#\{S_1[n_k,n_{k+1}]\cap\cdots\cap S_p[n_k,n_{k+1}]\}\ge\lambda \qquad \text{a.s.} \tag{5.7}$$

We consider the $2p$-dimensional random walk $\bar S(n)=(S_1(n),\dots,S_p(n))$. By the Markov property and Lévy's Borel–Cantelli lemma [see Corollary 5.29 in Breiman (1992)], (5.7) holds if we have

$$\sum_k \mathbb{P}^{\bar S(n_k)}\biggl\{J_{n_{k+1}-n_k}\ge\lambda\frac{n_{k+1}\log\log n_{k+1}}{(\log n_{k+1})^p}\biggr\}=\infty \qquad \text{a.s.} \tag{5.8}$$

Indeed, it is easy to see that $\sqrt{n_k \log\log n_k} = o(\sqrt{n_{k+1}/\log\log n_{k+1}})$ as $k \to \infty$. By the classic Hartman–Wintner law of the iterated logarithm, with probability 1 the events

$$\{\|\bar{S}(n_k)\| \le \sqrt{n_{k+1}/\log\log n_{k+1}}\}, \qquad k = 1, 2, \dots,$$

eventually hold. Therefore, (5.8) holds if we have

$$\sum_k \inf_{\|\bar{x}\| \le \sqrt{n_{k+1}/\log\log n_{k+1}}} \mathbb{P}^{\bar{x}}\left\{ J_{n_{k+1}-n_k} \ge \lambda \frac{n_{k+1}\log\log n_{k+1}}{(\log n_{k+1})^p} \right\} = \infty$$

which follows from Lemma 7 with $b_n = \log\log n$.

Since

$$J_{n_{k+1}} \ge \#\{S_1[n_k, n_{k+1}] \cap \cdots \cap S_p[n_k, n_{k+1}]\},$$

letting

$$\lambda \longrightarrow (2\pi)^p \left(\frac{2}{p}\right)^{p-1} \sqrt{\det(\Gamma)}\kappa(2,p)^{2p}$$

in (5.7) proves (5.1).

**Acknowledgments.** I thank R. F. Bass and J. Rosen for their interest in this work and for their helpful comments. I also wish to thank the referee for his valuable advice. He pointed out a mathematical error and an important reference.

Department of Mathematics
University of Tennessee
Knoxville, Tennessee 37996-1300
USA
e-mail: xchen@math.utk.edu
url: www.math.utk.edu/~xchen